\documentclass[12pt,reqno]{amsart}
\usepackage[margin=1in]{geometry}
\makeatletter
\def\section{\@startsection{section}{1}%
	\z@{.7\linespacing\@plus\linespacing}{.5\linespacing}%
	{\bfseries
		\centering
}}
\def\@secnumfont{\bfseries}
\makeatother
\usepackage{amsmath,amssymb,amsthm,graphicx,amsxtra, setspace}
\usepackage[utf8]{inputenc}
\usepackage{charter}
\usepackage{mathrsfs}
\usepackage{alltt}
\usepackage{stackengine}
\usepackage{relsize}
\usepackage{hyperref}
\usepackage{aliascnt}
\usepackage{tikz}
\usepackage{mathtools}
\usepackage{multicol}
\usepackage{upgreek}
\usepackage{graphicx,type1cm,eso-pic,color}
\usepackage{csquotes}
\allowdisplaybreaks
\usepackage{scalerel,stackengine}
\stackMath
\usepackage{setspace}
\def\intavg{\,\ThisStyle{\ensurestackMath{%
  \stackinset{c}{.2\LMpt}{c}{.2\LMpt}{\SavedStyle-}{\SavedStyle\phantom{\int}}}%
  \setbox0=\hbox{$\SavedStyle\int\,$}\kern-\wd0}\int}
\usepackage{scalerel}[2014/03/10]

\newtheorem{theorem}{Theorem}[section]
\newtheorem*{theorem*}{Theorem}

\newaliascnt{lemma}{theorem}
\newtheorem{lemma}[lemma]{Lemma}
\aliascntresetthe{lemma} 

\newaliascnt{proposition}{theorem}
\newtheorem{proposition}[proposition]{Proposition}
\aliascntresetthe{proposition}

\newaliascnt{assumption}{theorem}

\aliascntresetthe{assumption}

\newaliascnt{auxiliary}{theorem}

\aliascntresetthe{auxiliary}

\newaliascnt{corollary}{theorem}
\newtheorem{corollary}[corollary]{Corollary}
\aliascntresetthe{corollary}

\newaliascnt{definition}{theorem}

\aliascntresetthe{definition}

\newaliascnt{example}{theorem}

\aliascntresetthe{example}

\newaliascnt{remark}{theorem}
\newtheorem{remark}[remark]{Remark}
\aliascntresetthe{remark}

\newaliascnt{hypothesis}{theorem}

\aliascntresetthe{hypothesis}

\newaliascnt{property}{theorem}

\aliascntresetthe{property}

\DeclareMathOperator\supp{supp}
\usepackage[utf8]{inputenc}
\usepackage{amsmath}
\usepackage{enumerate}
\usepackage{bm}
\usepackage[T1]{fontenc}
\usepackage{amsthm}
\usepackage{amsfonts}
\usepackage{amssymb}
\usepackage{graphicx}
\usepackage{enumerate}
\usepackage[english]{babel}
\usepackage{comment}
\usepackage{amssymb}
\usepackage{amsthm}
\usepackage{xcolor}
\usepackage{dutchcal}
\usepackage{esint}

\usepackage[title]{appendix}
\usepackage{xcolor}

\DeclareMathOperator{\zmu}{z_{\mu}}
\def\vmu{{v_{\mu}}}
\def\vm{v_{\mu}}
\def\vmo{v_{\mu_0}}
\newcommand{\tvm}{\widetilde{v}_{\mu}}

\newcommand{\umr}{u_{\mu,R}}
\newcommand{\uM}{u_{M}}
\newcommand{\vM}{v_{M}}
\def\mpq{\mu^{p-1}+ \Tilde{C} \mu^{q-1}}
\def\Mpq{M^{p-1}+  M^{q-1}}

\DeclareMathOperator{\dus}{d^{s}_{U}}
\DeclareMathOperator{\durs}{d^{s}_{U_R}}
\def\kp{K_p}
\def\kq{K_q}
\def\kt{K_t}
\def\hp{H_p}
\def\hq{H_q}
\def\htt{H_t}

\def\vb{v_{\beta}}
\def\nn{\nu^{p-1}+  \nu^{q-1}}

\def\lo{\lambda^{p-1}+  \lambda^{q-1}}
\def\lp{\lambda^{p-1}+ (2C)^{\frac{q-p}{p-1}} \lambda^{q-1}}
\def\alo{ \alpha_1}
\DeclareMathOperator{\tail}{tail}
\DeclareMathOperator{\dis}{d}

\DeclareMathOperator{\daone}{d_{A_1}}
\DeclareMathOperator{\deone}{d_{E_1}}

\DeclareMathOperator{\dar}{d_{A_R}}
\DeclareMathOperator{\der}{d_{E_R}}
\DeclareMathOperator{\diss}{d^s}

\DeclareMathOperator{\deones}{d^{s}_{E_1}}

\DeclareMathOperator{\ders}{d^{s}_{E_R}}

\DeclareMathOperator{\doo}{d_\Omega}
\DeclareMathOperator{\doos}{d^{s}_{\Omega}}

\DeclareMathOperator{\dooso}{d^{-s}_{\Omega}}
\DeclareMathOperator{\osc}{osc}
\DeclareMathOperator{\loc}{loc}
\def\fp{(-\Delta)_{p}^{s}}
\def\fq{(-\Delta)_{q}^{s}}
\def\ft{(-\Delta)_{t}^{s}}
\def\brx{\Tilde{B}_{r,x_0}}
\def\bR{\Tilde{B}_{R}}
\def\bxR{\Tilde{B}_{x,R}}

\def\ex{ Ex(u,m,R)}
\def\exv{ Ex(v,m,R)}

\def\Ex{ Ex(u,M,R)}

\def\x{\Bar{x}}

\newcommand{\RR}{\mathbb{R}}
\newcommand{\Om} {\Omega}
\newcommand{\w}{\mathcal{W}}
\newcommand{\wo}{\mathcal{W}_0}
\newcommand{\tw}{\widetilde{\mathcal{W}}}
\newcommand{\two}{\widetilde{\mathcal{W}}_{0}}
\newcommand{\wsp}{W^{s,p}_0}
\newcommand{\wsq}{W^{s,q}_0}

\newcommand{\wpp}{W^{s,p}}

\newcommand{\va} {\varphi}
\newcommand{\var} {\varepsilon}

\newenvironment{sketch}{%
  \proof}{\endproof}
\newcommand{\Addresses}{{
		\footnote{
				\footnotesize
\noindent \textsuperscript{1}Indian Institute of Science Education and Research, Thiruvananthapuram 695551, India  \par\nopagebreak 
   \noindent \textsuperscript{2}Department of Mathematics, Indian Institute of Technology, Guwahati-781039, India  \par\nopagebreak 
\noindent \textsuperscript{A}\textit{e-mail:} \texttt{dhanya.tr@iisertvm.ac.in}.
\noindent \textsuperscript{B}\textit{e-mail:} \texttt{ritabrata20@iisertvm.ac.in}.
\noindent \textsuperscript{C}\textit{e-mail:} \texttt{uttam.maths@iitg.ac.in}.
\noindent \textsuperscript{D}\textit{e-mail:} \texttt{swetatiwari@iitg.ac.in}.

			\medskip\noindent
			{\bf Acknowledgments:} 
Ritabrata Jana was supported by the Prime Minister Research Fellowship during the execution of this research. R. Dhanya was supported by SERB MATRICS grant MTR/2022/000780 when this work was being carried out. We sincerely thank the anonymous reviewers for their insightful comments and suggestions, which have significantly enhanced the quality and clarity of this manuscript.  We would also like to extend our sincere gratitude to Professor Jacques Giacomoni of Université de Pau for his valuable discussions and insightful suggestions, which significantly contributed to enhancing the rigor of our work. His feedback was instrumental in refining the proofs and improving the overall presentation, particularly in the Appendix. 
			
}}}
\begin{document}
\title[Fine Boundary Regularity For The Fractional $(p,q)$-Laplacian]{ Fine Boundary Regularity For The Fractional $(p,q)$-Laplacian 	\Addresses	}
	\author[ R. Dhanya ]
	{ R. Dhanya\textsuperscript{1,A}} 
 	\author[Ritabrata Jana]
	{Ritabrata Jana\textsuperscript{1,B}} 
 \author[Uttam Kumar]
	{Uttam Kumar\textsuperscript{2,C}} 
 \author[Sweta Tiwari]
	{Sweta Tiwari\textsuperscript{2,D}} 
\maketitle
\begin{abstract}
In this article, we deal with the fine boundary regularity, a weighted H\"{o}lder regularity of weak solutions to the problem involving the fractional $(p,q)$ Laplacian denoted by
$(-\Delta)_{p}^{s} u + (-\Delta)_{q}^{s} u = f(x)$  in  $\Omega,$
and   $u=0$  in  $\mathbb{R}^N\setminus\Omega;$ 
where $\Omega$ is a $C^{1,1}$ bounded domain and $2 \leq p \leq q <\infty.$ For $0<s<1$ and for non-negative data $f\in L^{\infty}(\Omega),$  we employ the nonlocal analogue of the boundary Harnack method to establish that $u/{d_{\Omega}^{s}} \in C^{\alpha}(\Bar{\Omega})$ for some $\alpha \in (0,1),$ where $d_\Omega(x)$ is the distance of $x$ from the boundary. A novel barrier construction allows us to analyse the regularity theory even in the absence of the scaling or the homogeneity properties of the operator. Additionally, we extend our idea   to sign changing bounded $f$ as well and prove a fine boundary regularity for fractional $(p,q)$ Laplacian for some range of $s.$ 
\end{abstract}

\section{Introduction}

In recent years, significant attention has been devoted to research on nonlocal operators, with a notable emphasis on exploring its interior and boundary regularity results. Among these operators, particular attention has been devoted to  fractional $p$-Laplacian, defined as
\begin{equation*}
    \begin{aligned}
        (-\Delta)_p^s u(x):= 2 \underset{\varepsilon \rightarrow 0}{\lim} \int_{\RR^N \setminus B_{\varepsilon}(x)} \frac{|u(x)-u(y)|^{p-2}(u(x)-u(y))}{|x-y|^{N+ps}}\thinspace dy.
    \end{aligned}
\end{equation*}  
Such nonlocal operators find practical applications in real-world challenges, spanning domains like obstacle problems, finance, game theory, image processing, and so on. The Dirichlet problem associated with these operators is explored through the lenses of probability, potential theory, harmonic analysis and partial differential equations(see \cite{ Poz18} for further insights). In this article we 
establish an improved regularity result for the Dirichlet problem 
associated with a nonhomogeneous non-local operator known as the fractional $(p,q)$ Laplacian denoted by  $(-\Delta)_{p}^{s} u + (-\Delta)_{q}^{s} u .$

 Before presenting the main results of our article, we provide a concise overview of regularity results for nonlocal Dirichlet problems of the type $\mathcal{L}u = f $ in $\Omega,$ $u=0$ in $\Omega^c$  for various classes of $f.$ Within the realm of the nonlocal operators $\mathcal{L}$, the fractional p-Laplacian has been extensively studied.  When $p=2,$ the fractional p-Laplacian simplifies to the linear operator fractional Laplacian $(-\Delta)^s$ for which the regularity results are well understood. In the case when $f$ belongs to $L^\infty(\Omega)$, the weak solution of $(-\Delta)^s u=f$ in $\Omega$ under zero Dirichlet boundary condition belongs to $C^s(\mathbb{R}^n).$ This was proved by Ros-Oton and Serra in \cite{RS14} using the boundary Harnack method and the regularity is known to be optimal. Improved interior regularity results and Schauder-type estimates are also discussed in {\cite{{RS14}}} for the same Dirichlet boundary value problem under the assumption that $f\in C^\alpha(\Omega)$. Proceeding to nonlinear operators, Di Castro et al.in \cite{DKP16} and Cozzi \cite{Coz17}, explored the H\"{o}lder regularity and Harnack's inequality for minimizers of the equation $\mathcal{L} u =0$ within $\Omega$, where $u=g$ outside $\Omega$. Their focus was on nonlinear and homogeneous integro-differential operators, with the fractional p-Laplacian serving as the prototype model. Iannizzotto et al. \cite{IMS16} established the global H\"older continuity of weak solutions of $(-\Delta)_p^s u = f$ in $\Omega$, $u=0$ in $\Omega^c$, specifying that $u\in C^\beta(\overline{\Omega})$ for some $\beta\in (0,1).$  This result was demonstrated for bounded measurable functions $f$, although the H\"older exponent $\beta$ was unspecified in \cite{IMS16}. Later, Iannizzotto et al.\cite{IMS20} remarked that the solution $u$ indeed belongs to $C^s(\mathbb{R}^n)$ itself, implying that the exponent $\beta$ may be chosen as $s$ (refer to Theorem 2.7 of \cite{IM24} for a proof).

Researchers have also established various regularity results for weak solutions of the equation $(-\Delta)_p^s u=f$, where the given data $f$ is not necessarily a bounded function but satisfies certain integrability conditions. In \cite{BP16}, Brasco et.al. proved $L^\infty$ regularity and continuity of weak solutions when $f$ belongs to $L^q(\Omega)$ for large exponents $q.$ Furthermore, in \cite{BLS18}, authors proved stronger results such as the optimal interior H\"older continuity of weak solutions for the Dirichlet problem involving the fractional p-Laplacian when $p\geq 2$ and in \cite{GL23} for the case $1<p<2.$ Higher Sobolev and H\"older regularity results are proved under various integrability conditions in \cite{BL17}. Furthermore, in their work, Kuusi et. al. \cite{KMS15} explored the existence and regularity results for the solutions of the equations modelling fractional $p$-Laplace problems with measure data.

Regularity results are also explored for more general class of nonlocal linear operators with the kernels which are not translation invariant. Cafarelli and Silvestre achieved significant progress in understanding interior regularity, as demonstrated in \cite{CS11} where they obtained $C^{1,\alpha}$ interior regularity for the viscosity solutions through approximation techniques. For a broader class of linear fractional problems involving more general, possibly singular, measurable kernels, Dyda and Kassmann's work \cite{DB20} provides the corresponding regularity results. Bonder et al. \cite{BFS22} have established a global H\"{o}lder continuity result for weak solutions in problems featuring the fractional $(-\Delta)^{s}_{g}$ Laplacian, where $g$ represents a convex Young's function. Readers are also referred to \cite{Fal20,Kas09,Mos18,Sil06}, for  the further insights into H\"older regularity and to \cite{Bas09,FR23} for the Schauder estimates for nonlocal operators with more general kernels.  

Aforementioned discussions on the regularity have largely focused on nonlocal operators which observe homogenity and scaling properties. 
However, our interest lies in the fractional $(p,q)$ Laplace operator, which does not exhibit these characteristics.  This operator represents a specific category of nonuniformly elliptic problems and are relevant in various practical applications, such as in the homogenization of strongly anisotropic materials. In \cite{DP19}, the authors have successfully derived interior H\"{o}lder regularity results for viscosity solutions within a more general class of fractional double-phase problems.  For fractional operators with non-standard growth and with homogeneous data, local boundedness and H\"older continuity have been proved in \cite{BJS22,YZ23}. Recent studies have delved into the interior and boundary regularity of weak solutions to the fractional $(p,q)$ Laplace operator with Dirichlet boundary condition, as noted in articles \cite{GKS20,GKS22,GKS23}. To be precise, Giacomoni et. al. \cite{GKS22} considered the problem 
 \begin{equation}
\label{neweqn111}
 \left\{
     \begin{array}{rll}
         (-\Delta)^{s_1}_p  u + (-\Delta)^{s_2}_q u \thinspace & = &\thinspace f \text{ in } \Omega,
         \\
         u \thinspace& =& \thinspace 0   \text{ in } \Omega^{c},
     \end{array}\right.
 \end{equation}
 where $2\leq p \leq q <\infty, 0\leq s_1\leq s_2\leq 1$ and $\Om$ a bounded domain in $\RR^N$ with $C^{1,1}$ boundary. In Remark 12 of \cite{GKS22}, authors proved that the solution to \eqref{neweqn111} lies in the function space $C^s(\mathbb{R}^N)$ when $s_1=s_2=s.$ Furthermore, if $f\in L^\infty(\Omega)$ from  Remark 3.7 of  \cite{GKS23}, it is evident that $\frac{u}{\doos}$ belongs to $L^\infty(\Omega)$. Our objective in this article is to demonstrate the H\"{o}lder regularity of $\frac{u}{\doos}$  where $u$ is the weak solution of \eqref{neweqn111} when $s_1=s_2=s.$  
\begin{theorem}\label{mainthmst}
     Let $\Omega$ be a bounded domain in $\mathbb{R}^N$ with $C^{1,1}$ boundary and $\mathrm{d}_{\Omega}(x)$ is defined in \eqref{disdef}.
     Let $u\in W^{s,p}_0(\Omega)\cap W^{s,q}_0(\Omega)$ be the weak solution of 
     \begin{equation}
\label{eqn1}
\left\{
    \begin{array}{rll}
        \fp u + \fq u \thinspace & = & \thinspace f \text{ in } \Omega,
        \\
       u \thinspace & = & \thinspace 0   \text{ on } \Omega^{c},
    \end{array}\right.
\end{equation}
where $2 \leqslant p \leqslant q<\infty.$ Assume further that any of the following two conditions are satisfied:
 \begin{itemize}
         \item [(a)] $f\in L^\infty(\Omega),$ $f\geq 0$ and $0<s<1,$
         \item [(b)] $f\in L^\infty(\Omega)$ and $0<s<\frac{1}{q}.$
     \end{itemize}
     Then, $\frac{u}{\doos}\in C^\alpha(\bar{\Omega}) $ for some $\alpha\in (0,s]$ and $\displaystyle \left\|\frac{u}{\doos}\right\|_{C^\alpha(\bar{\Omega})} \leq C$ where $C$ depends on $N, p,q, s$, $\Omega$ and $\|f\|_{L^{\infty}(\Omega)}.$ 
\end{theorem}

 The significance of the fine boundary regularity result for nonlocal problems lies in its pivotal role in establishing $C_{s}^{\alpha}(\overline{\Omega})$ as the appropriate function space for the solution of many elliptic and parabolic problems
where
 \begin{equation*}
     \begin{aligned}
 C_{s}^{\alpha}(\Bar{\Omega}):=\left\{u \in C^{0}(\Bar{\Omega}): \frac{u}{\diss} \text{ has a } \alpha-\text{H\"{o}lder continuous extension to }\Bar{\Omega}\right\}.        
     \end{aligned}
 \end{equation*}
In the local case, the corresponding regularity results for the $p$-Laplacian have been explored by Lieberman in \cite{Lib88,Lib91}, where the counterparts of these $C_{s}^{\alpha}$ spaces are the $C^{1,\alpha}$ spaces. Various applications of $C^{1,\alpha}$ regularity result for quasilinear problems are now standard, with some of these results documented in \cite{MZ97}. As an application of regularity results for nonlocal problems, exploiting the behaviour of $u/d^s$ near the boundary, Ros-Oton et al.\cite{RS14b} have proved the Pohozaev identity for the fractional-Laplacian. Additionally, this improved  boundary regularity result is utilized in \cite{Aba13} to ensure the well-definedness of test functions and to investigate blow-up solutions for the fractional Laplacian.  Fall et al. in \cite{FJ15} leverage the fine boundary regularity  to explore overdetermined problems in the nonlocal settings. In the context of nonlinear problems involving fractional operators, one can leverage the benefits of the compact embedding of $C_{s}^{\alpha} \hookrightarrow C_{s}^{0}$ to prove various existence and multiplicity results. The existence of solutions for the fractional $p$-Laplacian using a variational approach is demonstrated in \cite{IMS20b, IMP23}, while discussions in the realm of degree theory are presented in \cite{FI23}.  By employing the idea of compact embedding, the multiplicity of solutions can also be established, as evidenced in \cite{IL21, FI21, FI22}. We encourage the readers to refer \cite{IM24} for more applications of fine boundary regularity results concerning the fractional $p$-Laplacian. In our previous work \cite{DJKT23}, we investigated the existence of positive solutions for semipositone problems involving the fractional p-Laplacian, relying heavily on the fine boundary regularity results of Iannizzotto \cite{IMS20}. Upon establishing the main result of this paper, we plan to extend this problem to incorporate the fractional $(p,q)$ Laplace operator. We believe that the crucial ideas involved in the proof of this article can be applied to establish the fine boundary regularity result for $(p, q)$ Laplacian  when both $p,q$ are less than 2, following the approach of \cite{IM24}. This will be addressed in a subsequent work.

We shall now briefly discuss the idea behind the proof of the improved boundary regularity result  for nonlocal problems.
 Ros-oton et al.\cite{RS14} have proved the fine boundary regularity result for the fractional Laplacian (i.e. $p=q=2$) by \enquote{trapping the solution}  between two multiples of $\diss$ in order to control the oscillation of $\frac{u}{\diss}$ near the boundary. They have proved a fractional analog of boundary Harnack inequalities by constructing proper upper and lower barrier for the problem.  Later, Iannizzotto et al. \cite{IMS20} have followed a similar approach to extend the result for fractional $p-$ Laplacian for $p \geq 2,$ where two weak Harnack inequalities have been established for the function $\frac{u}{\doos}$: one in the case where $u$ is a subsolution, and another in the case where $u$ is a supersolution. The study aims to control the behaviour of $\frac{u}{\doos}$ near the boundary through the nonlocal excess, defined as 
 \begin{equation*}
     \begin{aligned}
         \operatorname{Ex}(u,m,R,x_0)= \intavg_{\brx} \left|\frac{u(x)}{\doos(x)}-m\right|\thinspace dx,
     \end{aligned}
 \end{equation*}
 where $x_0\in \partial \Omega,$ $ m\in \RR$, $R>0.$ Here, $\brx$ represents a small ball with a radius comparable to $R$, situated at a distance greater than $R$ in the normal direction from $x_0\in \partial\Om$, and the ball is positioned away from the boundary (see Figure 1). The term \enquote{nonlocal excess} is used because, given a bound on $\fp u+\fq u $, the pointwise behavior of $\frac{u}{\doos}$ inside $ B_{R}(x_0)\cap \Omega$ is determined by the magnitude of the excess of $u$ in $\brx.$ 
 \par
 Due to the lack of additivity properties of the fractional p-Laplacian, Iannizzotto et al. extensively focuses on constructing two families of one-parameter basic barriers $w_\lambda,$ which are defined in (1.6) and pp. 6 of \cite{IMS20} where $\lambda\simeq Ex(u).$

For small values of $\lambda,$ they have constructed the barrier starting from $\doos$ and performing a $C^{1,1}$ small diffeomorphism. But for the large values of $\lambda, $ they have taken advantage of the homogeneity and scaling properties of the operator and yields the barrier as a multiple of the torsion function which is defined as the unique solution to the problem 
\begin{equation*}\left\{
    \begin{aligned}
        &\fp v = 1 \text{ in }  B_{R/2}(x)\cap \Omega,\\
        & v=0  \text{ in } (B_{R/2}(x)\cap \Omega)^{c}.
    \end{aligned}\right.
\end{equation*}
Clearly, a multiple of the torsion function satisfies, $(-\Delta)_p^s (\lambda^{\frac{1}{p-1}} v)=\lambda$ in $ B_{R/2}(x)\cap \Omega$ thus simplifying the analysis of its behavior, and, consequently, the construction of barriers. Notably, when dealing with the fractional $(p,q)$ Laplacian operator, we encounter a limitation in adopting the same approach as in previous studies.
The novelty of our work lies in constructing such barrier when the operator lack homogenity and scale invariance property. We overcome this difficulty by estimating the asymptotic behavior of solutions to a family of nonlocal PDE's depending on a  parameter $\lambda$, utilizing the regularity results established by Giacomoni et al. in \cite{GKS20} and \cite{GKS23}.

\par The article is structured as follows: Section \ref{prelim} provides definitions and notations essential for understanding the subsequent sections. In Section \ref{est}, we establish the asymptotic behaviour of solutions mentioned in the above paragraph. Following this, Section \ref{lowbdsec} presents a lower bound for supersolutions of the fractional $(p,q)$ Laplacian, while Section \ref{upbdsec} focuses on proving an upper bound for subsolutions when $f\geq 0$. Subsequently, Section \ref{whrp} is dedicated to demonstrating an oscillation bound and weighted H\"older regularity result for solutions of the fractional $(p,q)$ Laplacian with non-negative data. In Section \ref{whrs}, we establish the fine boundary regularity result for sign-changing bounded data when $s\in(0,\frac{1}{q})$. Lastly, Appendix \ref{app} we present the  auxiliary results needed for the paper.
\section{Preliminaries}\label{prelim}
 We begin this section by introducing definitions which are relevant for this article. For a measurable function $u:\mathbb{R}^{N}\rightarrow \mathbb{R}$, we define Gagliardo seminorm
\begin{equation*}
	[u]_{s,t}:=[u]_{W^{s,t}(\mathbb{R}^N)}:=  \left(\int_{\mathbb{R}^N \times \mathbb{R}^N} \dfrac{|u(x)-u(y)|^{t}}{|x-y|^{N+st}}\,dx\, dy\right)^{1/t},
\end{equation*}
for $1<t<\infty$ and $0<s<1.$
We consider the space $W^{s,t}(\mathbb{R}^N)$ defined as  
\begin{equation*}
	W^{s,t}(\mathbb{R}^{N}):= \left \{u \in {L}^{t}(\mathbb{R}^{N}):[u]_{s,t}<\infty  \right\}.
\end{equation*}
The space $W^{s,t}(\mathbb{R}^{N})$ is a Banach space with respect to the norm 
\begin{equation*}
\|u\|_{{W^{s,t}(\mathbb{R}^{N})}}=\left( \|u\|^{t}_{L^{t}(\mathbb{R}^{N})} + [u]^{t}_{{W^{s,t}(\mathbb{R}^{N})}}\right)^{\frac{1}{t}} . 	
\end{equation*}
A comprehensive examination of the fractional Sobolev Space and its properties are presented in \cite{DPV12}.
Let $\Om \subset \RR^N$ be a bounded domain with a $C^{1,1}$ boundary. To address the Dirichlet boundary condition, we naturally consider the space $W^{s,t}_0(\Omega)$ defined as
\begin{equation*}
	W_{0}^{s,t}(\Omega):= \left \{u \in W^{s,t}(\mathbb{R}^{N}):u=0\medspace\text{in}\medspace \mathbb{R}^{N}\setminus \Omega \right\}.
\end{equation*}
This is a separable, uniformly convex Banach space endowed with the norm $$\|u\|=	\|u\|_{{W^{s,t}(\mathbb{R}^{N})}}.$$ Moreover the embedding $W^{s,t}_{0}(\Omega)\hookrightarrow L^{r}(\Omega)$ is continuous for $1\leq r\leq t^{*}_{s}:=\frac{Nt}{N-ts}$ and
compact for  $1\leq r < t^{*}_{s}$. 
 Due to continuous embedding of $W^{s,t}_{0}(\Omega)\hookrightarrow L^{r}(\Omega)$ for $1\leq r\leq t^{*}_{s}$, we define the equivalent norm on $W^{s,t}_{0}(\Omega)$ as
\begin{equation*}
	\|u\|_{W^{s,t}_{0}}:=\left(\int_{\mathbb{R}^N \times \mathbb{R}^N} \dfrac{|u(x)-u(y)|^{t}}{|x-y|^{N+st}}\,dx\, dy\right)^{1/t} .
\end{equation*} 
The dual space of $W^{s,t}_{0}(\Omega)$ is denoted by $W^{-s,t'}(\Om)$ for $1<t<\infty.$
We shall also recall the following space from \cite{IMS16}
\begin{equation*}
\begin{aligned}
    \widetilde{W}^{s,t}(\Om):= \Bigg\{u \in L_{\loc}^{t}(\RR^N): \exists\; \Om^{\prime} \supset\supset \Om &\text{ such that } u\in W^{s,t}(\Om^{\prime})
   \\ &\text{ and } \int_{\RR^N}\frac{|u(x)|^{t-1}}{(1+|x|)^{N+ps}}<\infty \Bigg\}. 
\end{aligned}            
\end{equation*}
Using \cite[Lemma 2.3]{IMS16} one can prove that if $u \in  \widetilde{W}^{s,t}(\Om),$ then $(-\Delta)_{t}^{s}u \in W^{-s,t'}(\Om).$ Let $U$ be an open subset of $\Om$ and we set 
\begin{equation*}
    \begin{aligned}
        \widetilde{W}_{0}^{s,t}(U)= \left\{u \in \widetilde{W}^{s,t}(U): u=0 \text{ in } \Om^c  \right\} .
    \end{aligned}
\end{equation*}
In \cite{GKS20}, it is proved that $ W^{s_2,q}(\Omega)\hookrightarrow W^{s_1,p}(\Omega)$ for $1<p \leq q < \infty,$  $0<s_1< s_2<1$ and when $\Omega$ is bounded. However, the result is not true for the case $s_1=s_2,$ a counter example is provided in \cite[Theorem 1.1]{MS15}. Because of this lack of continuous embedding of the space $\wsp$ and $\wsq$, in order to consider the weak solution associated with the operator $(-\Delta)_p^s+ (-\Delta)_q^s$  we consider
\begin{equation*}
    \mathcal{W}_{0}(\Om)=  W_{0}^{s,p}(\Om)\cap W_{0}^{s,q}(\Om) ,
\end{equation*} 
with the norm $\|\cdot\|_{\mathcal{W}_{0}(\Om)}=\|\cdot\|_{W_{0}^{s,p}(\Om)}+\|\cdot\|_{W_{0}^{s,q}(\Om)}$. The space $\mathcal{W}^{\prime}(\Om)$ is denoted for the dual of the $\mathcal{W}_{0}(\Om).$ Analogously we define 
\begin{equation*}
\begin{aligned}
    &\widetilde{\mathcal{W}}(\Om)=  
     \widetilde{W}^{s,p}(\Om)
     \cap\widetilde{ W}^{s,q}(\Om),
      \\
    &\widetilde{\mathcal{W}}_{0}(\Om)=  
     \widetilde{W}_{0}^{s,p}(\Om)
     \cap\widetilde{ W}_{0}^{s,q}(\Om).
\end{aligned}
\end{equation*}
Drawing the inspiration from \cite{IMS20}, for all $R>0$ and $t\geq 1,$ we define the  nonlocal tail for any measurable function $u: \RR^N \mapsto \RR$ as 
\begin{equation*}
    \begin{aligned}
        \tail_{t}(u,R)=\left[\int_{\Om\cap B_{R}^{c}(0)}\frac{|u(x)|^{t}}{|x|^{N+s}}\thinspace dx\right]^{\frac{1}{t}} .
    \end{aligned}
\end{equation*}
Note that the tail defined here is different from the notion of tail defined in \cite{GKS23} by a factor of $R^{s}.$ Next, we provide a set of notations concerning specific subsets of $\Omega$ which are used in the subsequent analysis of the paper. For all $x\in \RR^N$ and $R>0$ we set $D_R(x):= B_R(x)\cap\Om$ where $B_R(x)$ is a ball of radius $R$ centered at $x.$ When the centre is the origin, we may denote it by $B_R$ and $D_R.$ The distance function $\dis_{U}: \RR^N \rightarrow \RR_{+}$ is defined as 
\begin{equation}\label{disdef}
    \dis_{U}(x):=\underset{y \in U^c}{\inf \thinspace } |x-y| \text{ where } U \text{ is an open set in } \RR^{N}.
\end{equation}

 As we assume the boundary $\partial\Omega$ has $C^{1,1}$ regularity, $\doo$ is a Lipschitz continuous function in $\mathbb{R}^N.$ Moreover the interior sphere property holds for the domain $\Omega.$ Specifically, there exists $R > 0$ such that for each $x \in \partial \Omega$, there exists $y \in \Omega$ such that the ball $B_{2R}(y)$, tangent to $\partial \Omega$ at $x$, is contained entirely within $\Omega$. Define
\begin{equation*}
    \begin{aligned}
        \rho:=\rho(\Om):=\sup\left\{R \thinspace: \forall \ x \in \partial\Omega \ \exists \ B_{2R}(y)\subseteq \Om \text{ s.t. } x\in \partial B_{2R}(y)\right\}.
    \end{aligned}
\end{equation*}
This $\rho$ represents the supremum of all such radii $R$ for which the mentioned tangential ball inclusion property holds for points on the boundary of $\Omega$. Next we define $\rho$-neighborhood of $\partial\Omega$ as
\begin{equation*}
    \begin{aligned}
        \Om_{\rho}:=\left\{x \in \Om \ : \ \doo(x)< \rho\right\}.
    \end{aligned}
\end{equation*}
\begin{center}
\includegraphics[scale = 0.4]{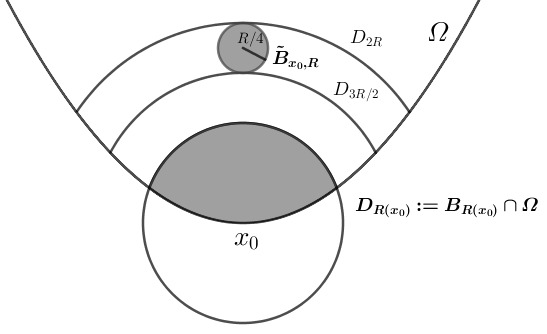}    
\end{center}
With the chosen value of $\rho$, the metric projection $\Pi_{\Omega}:\Omega_{\rho}\rightarrow \partial\Omega$ is well-defined and constitutes a $C^{1,1}$ map. Moreover, for all $x\in \partial\Omega$ and $R \in (0,\rho)$ there exists a ball $\bxR$ of radius $R/4,$ s.t.
\begin{equation*}
    \begin{aligned}
        \bxR\subset D_{2R}(x)\setminus D_{3R/2}(x), \thickspace \underset{ y \in \bxR}{\inf\thinspace} \doo(y) \geqslant 3R/2.
    \end{aligned}
\end{equation*}
If $x=0,$ we denote $\bxR$ as $\bR.$ We define our nonlocal excess as 
\begin{equation*}
    \operatorname{Ex}(u,m,R,x_0)= \intavg_{\brx} \left|\frac{u(x)}{\doos(x)}-m\right|\thinspace dx.
\end{equation*} 
Naturally $\operatorname{Ex}(u,m,R,0)$ is denoted as $\ex.$
\par 
A function $u\in \wo(\Om)$ is a weak supersolution (sub olution) of \eqref{eqn1} if	
	\begin{equation*}
		\sum_{t=p,q}\; \int_{\mathbb{R}^N \times \mathbb{R}^N } \dfrac{|u(x)-u(y)|^{t-2}(u(x)-u(y))(\varphi(x)-\varphi(y))}{|x-y|^{N+st}}  \,dx\,dy 
		 \geq (\leq) 	\int_{\Omega} f(x) \varphi(x) \,dx,	
	\end{equation*}
for every $\varphi \in 	\wo(\Omega)_{+}$.	
 \par Throughout the paper, all equations and inequalities involving $\fp+ \fq $ are understood in the weak sense. In this manuscript, we adhere to the methodologies established by Ros-Oton et al. \cite{RS14} and Iannizzotto et al. \cite{IMS20}, with the exception of relying on the homogeneity and scaling properties of the operator. While we maintain alignment with their fundamental principles and cite their work when relevant, it is crucial to emphasize that our operator is nonhomogeneous and hence we present the detailed calculations.

\underline{\textbf{Notation}}
\\
Throughout this paper $\Om$ is a bounded domain with $C^{1,1}$ smooth boundary and $2 \leqslant p \leqslant q < \infty.$ In this article, the real parameter $s$ lies in the interval $(0,1)$, except for section 7, where $s\in (0,\frac{1}{q})$. Unless stated otherwise, $k,M$ and $C$ etc. represent generic positive constants. For $S\subset\RR^{2N},$ we denote
\begin{equation}
    \begin{aligned}
        A_{t}(u,v,S):=\int_{S}\frac{|u(x)-u(y)|^{t-2}(u(x)-u(y))(v(x)-v(y))}{|x-y|^{N+st}} \thinspace dx\thinspace dy.
    \end{aligned}
\end{equation}

\section{Some Important Estimates}\label{est} 
{
In this section, we will present two lemmas that illustrate the precise behavior of a family of solutions for a class of boundary value problems involving the $(p,q)$-Laplacian. As mentioned in the introduction, to prove the fine boundary regularity result, we aim to control the oscillation of $\frac{u}{\doos}$ near the boundary, when $u$ is either a sub or supersolution of \eqref{eqn1}. Previously, for fractional p Laplace operator, the torsion function played an important role in controlling the oscillations of $\frac{u}{\doos}.$ Now, due to the absence of scaling and homogeneity properties of the $(p,q)$-Laplacian  the direct application of the torsion function is not feasible. Instead, the estimates we derive in this section  will offer alternative pathways to obtain the desired control.
}
\par
Let $U\subset \Omega $ be a bounded open set with $C^{1,1}$ boundary. Define $U_R:=\{x \,| \, \frac{x}{R}\in U\}.$
Then, by definition, we have
\begin{equation}\label{diseqn1}
    \begin{aligned}
        R \dis_{U}(x) = \dis_{U_R} (Rx) \text{ for } x\in U.
    \end{aligned}
\end{equation}
In the forthcoming sections, we will use the results we prove here with the domain $U$ replaced by $E, A$ or $D$ as applicable.

Firstly, we state a lemma which establishes that the solution to equation \eqref{eextlem} within the domain $U$ acts as a subsolution across the entire domain $\mathbb{R}^N$. The proof relies on the  convexity of the energy functional associated with the fractional $(p,q)$-Laplacian operator. We omit a detailed proof, as it follows a similar
 approach outlined in \cite[Lemma 2.4]{IMS20}.
\begin{lemma}\label{extlem}
    Let  $C>0$,  and $v \in \wo(U)$ solves 
    \begin{eqnarray}\label{eextlem}
    \left\{
	\begin{array}{rll}
	\fp v + \fq v &= C &\text{ in } U,
            \\
            v&=0 &\text{ in } U^c.
	\end{array} 
	\right. 
\end{eqnarray}
    Then $\fp v + \fq v \leq C$ in $\RR^N.$ 
\end{lemma}
\begin{lemma} \label{lowbdlem}
    Let $\umr \in \wo(U_R)$  be the solution of
        \begin{eqnarray}\label{umreqn}
    \left\{
	\begin{array}{rll}
  (-\Delta)_{p}^s \umr(y)+ (-\Delta)_{q}^s \umr(y) &=\frac{\mpq}{R^s}  &\text{for} \ y \in U_R,\\
\umr(y) &= 0  &\text{for} \ y\in U_{R}^{c},
	\end{array} 
	\right. 
\end{eqnarray}
where $0< R < \rho/4,$ and $\Tilde{C}$ be a given  positive constant. For any $ \bar{\mu}_{0}>0,$  there exists $C_{1}(\bar{\mu}_{0},\Tilde{C},U)>0$ such that if $\mu>\bar{\mu}_{0}$ then  
 \begin{equation*}
      \frac{\mu}{C_1} \durs(y) \leq \umr(y) \leq C_1 \mu R^s \text{ for } y \in U_R,
 \end{equation*} 
where the constant $C_1$ is independent of $\mu, R.$ 
\end{lemma}
\begin{proof}
 We first claim that there exists a $ \mu_{0}>0$ and a constant $C_{1}(\mu_{0},\Tilde{C},U)>0$ such that if $\mu>\mu_{0}$ then  
 \begin{equation}\label{claimmu}
      \frac{\mu}{C_1} \durs(y) \leq \umr(y) \leq C_1 \mu R^s \text{ for } y \in U_R.
 \end{equation}       We first set $v_{\mu_1}(x)=\umr(Rx)$ for $x\in U.$ Then, 
 \begin{equation*}
     \begin{aligned}
            R^{(q-p)s} \fp v_{\mu_1}(x) + \fq v_{\mu_1}(x) =  R^{(q-1)s} (\mpq) \text{ for } x\in U. 
     \end{aligned}
 \end{equation*}
 Next we define $\tvm(x)= \frac{v_{\mu_1}(x)}{\mu R^s},$ for $x \in U.$ Then, $\tvm$ satisfies 
\begin{eqnarray}\label{veqn}
    \left\{
    \begin{array}{rll}
    \mu^{p-q} \fp \tvm + \fq \tvm &= \mu^{p-q}+ \Tilde{C} &\text{ in }U,
         \\
         \tvm&=0 &\text{ in }U^{c}.
    \end{array}
    \right.
\end{eqnarray}
We claim that for any $K_1\Subset U,$ $\inf_{K_1} \tvm\geq C(\mu_0)>0$ independent of $\mu$ and $R.$ To prove the claim, we use $\tvm$ as a test function in \eqref{veqn}, we have 
    \begin{equation*}
         \relax  [\tvm]_{\wsq}^q \leq  \mu^{p-q} [\tvm]_{\wsp}^p+[\tvm]_{\wsq}^q\leq  (\Tilde{C}+\overline{\mu}_0^{p-q}) \int  \tvm.
    \end{equation*}
    Thanks to the continuous embedding $\wsq(U) \hookrightarrow L^1(U)$ for any $q>1,$ we get $[\tvm]_{W_{0}^{s,q}}\leq C,$ independent of $\mu$ for large $\mu.$  {Hence upto a subsequence, $\tvm\rightharpoonup v_0$ in $\wsq(\Om)$. Using Theorem \ref{csreg} in Appendix, we have that $\|\tvm\|_{C_{\loc}^{s}}$ is uniformly bounded for large values of $\mu\geq \mu_{0}>0.$ Hence, we can show that 
    \begin{equation*}
        \begin{aligned}
           \left| \int_{\mathbb{R}^N \times \mathbb{R}^N } \dfrac{(\tvm(x)-\tvm(y))^{p-1}(\varphi(x)-\varphi(y))}{|x-y|^{N+sp}}  \,dx\,dy \right| \leq C 
        \end{aligned}
    \end{equation*}
    for each $\varphi \in C_c^\infty(\Om)$ where $C(\mu_{0})$ is independent of $\mu.$ This clearly implies that as $\mu\rightarrow \infty,$
     \begin{equation*}
        \begin{aligned}
            \mu^{p-q}\int_{\mathbb{R}^N \times \mathbb{R}^N } \dfrac{(\tvm(x)-\tvm(y))^{p-1}(\varphi(x)-\varphi(y))}{|x-y|^{N+sp}}  \,dx\,dy \rightarrow 0 \mbox{ for all } \varphi\in C_c^\infty(\Omega).
        \end{aligned}
    \end{equation*}
    Due to the weak-weak continuity property \cite[Lemma 2.2]{CMM18} of fractional $q$-Laplacian and the density argument, letting $\mu\rightarrow \infty$ in (\ref{veqn}) we can prove that the function $v_0$ solves
    \begin{equation}\label{evo}
        \begin{aligned}
            \fq v_0 &=\Tilde{ C} \ \text{in} \ U,\\
v_0 &= 0 \hspace{0.2cm} \text{in} \ U^{c}.
        \end{aligned}
    \end{equation}
   Thanks to the $C_{\loc}^{s}(U)$ uniform bound(Theorem \ref{csreg}), we can apply the Ascoli-Arzelà theorem to deduce that $\tvm$ has a uniformly convergent subsequence on every compact subset of $U.$ Uniqueness of the solution of \eqref{evo} would imply that   $\|\tvm - v_0\|_{L^\infty(K_1)}\rightarrow0$  as $\mu \to \infty$ for any fix compact subset $K_1$ of $U.$  
\\
Now, using the strong comparison principle for fractional $q$-Laplacian \cite[Lemma 2.3]{IMS20},  we have  $\underset{K_1}{\inf} \, v_0\geq k (\Tilde{C})^{\frac{1}{p-1}}= C_1 >0$ for the same $\Tilde{C}$ given in \eqref{evo} and some constant $k(N,s,q,K_{1})>0.$ Since $\tvm$ uniformly converges to $v_0,$  there exists a $\mu_0>0$ such that 
   \begin{equation*}
    \begin{aligned}
      \tvm(x)\geq v_0(x)-\varepsilon \geq C_1-\varepsilon \ \text{for all } \ x\in  K_1, \text{ and for all }\mu\geq \mu_0.
    \end{aligned}
\end{equation*}
 Now, choosing $\varepsilon$ small enough we obtain, 
\begin{equation}\label{Ceqn}
    \begin{aligned}
        \underset{K_1}{\inf} \, \tvm \geq C(\mu_{0})>0  \, \text{  for all } \mu\geq \mu_0 \,\, \text{ and } \forall x\in K_1 .
    \end{aligned}
\end{equation}}
 Next to derive an estimate near the boundary of $U$, we set $w(x)= k \dus(x)$ for $x \in U$ and $0<k\leq \frac{C}{\underset{\RR^N}{\max} \dus}$ where $C$ is same as given in \eqref{Ceqn}. We know that $$\underset{U\setminus K_1}{\max} \{|\fp \dus|, |\fq \dus| \} \leq M_1$$ for some constant $M_1>0.$ If we choose $k$ such that $\mu_{0}^{p-q}k^{p-1}M_1 + k^{q-1}M_1 \leq \Tilde{ C}/2,$ along with the condition $0<k\leq \frac{C}{\underset{\RR^N}{\max} \dus}$,   then 
 \begin{eqnarray*}
     \begin{array}{rll}
          \mu^{p-q} \fp (k \dus(x)) + \fq (k \dus(x)) & \leq \mu^{p-q}+ \Tilde{C} &\text{ in } U\setminus K_1,
         \\
        k \dus(x) \leq C & \leq \underset{K_1}{\inf}\, \widetilde{v}_{\mu} \leq \widetilde{v}_{\mu} &\text{ in } K_1. 
     \end{array}
 \end{eqnarray*}
By applying the comparison principle, we obtain $\widetilde{v}_{\mu} \geq k \dus$ in $\RR^N,$ implying that 
 \begin{equation*}
     \umr(Rx) \geq k \mu R^s \dus(x), \text{ for } x \in U.
 \end{equation*}
 Combining with \eqref{diseqn1} choosing $C_1$ large enough, we get that
 \begin{equation}\label{lumr}
     \begin{aligned}
         \umr(y)\geq \frac{\mu}{C_1} \durs(y) \text{ for } y \in U_R \text{ and } \mu\geq \mu_0.
     \end{aligned}
 \end{equation}
Now we shall find the upper bound for $u_{\mu, R}.$ Let $\vmo$ be the solution of the following equation: 
    \begin{eqnarray*}
        \left\{
        \begin{array}{rll}
          \mu^{p-q} \fp \vmo(x) + \fq \vmo(x) &=    \mu_{0}^{p-q}+\Tilde{C} \ &\text{in} \  U,\\
\vmo &= 0 \ &\text{on} \ U^{c}.
        \end{array}
        \right.
    \end{eqnarray*}
    Thanks to \cite[Theorem 2.3]{GKS20}, we have $\|\vmo\|_{L^\infty(U)} \leq M$ where $M$ is independent of $\mu.$ Since,
\begin{equation*}
    \mu_0^{p-q}+\Tilde{C} >\mu^{p-q} +\Tilde{C},
\end{equation*} 
using the comparison principle for $\tvm$ and $\vmo,$  we have $\underset{U}{\sup \thinspace } \tvm\leq M.$  Since $\tvm(x)=\frac{\umr(Rx)}{\mu R^s}$ for all $x\in U$ we conclude that 
\begin{equation}\label{euumr}
    \begin{aligned}
        \umr(y)\leq M \mu R^s \text{ in } U_R \text{ and } \mu\geq \mu_0.
    \end{aligned}
\end{equation}
The claim \eqref{claimmu} now immediately follows from \eqref{lumr} and \eqref{euumr}.
\\
For any fixed \(\bar{\mu}_{0} > 0\), if \(\bar{\mu}_{0} < \mu <\mu_0\),  by applying the weak comparison principle, we obtain \( u_{\bar{\mu}_0,R}\leq u_{\mu, R} \leq u_{{\mu}_{0}, R} \). This observation along with the main claim \eqref{claimmu} allow us to conclude the desired result.
\end{proof}
Now we remove the $\mu_{0}$ dependence of the constant, but this comes at the cost of the result being valid only within a specific range of $s$.
\begin{lemma}\label{muuniformlemma}
Let $0<s<\frac{1}{q}$ and  $\umr \in \wo(U_R)$  solves \eqref{umreqn}
 for $0< R < \rho/4,$  $\mu>0$ and $\Tilde{C}$ be a given  positive constant.
 Then 
 \begin{equation*}
      \frac{\mu}{C_1} \durs(y) \leq \umr(y) \leq C_1 \mu R^s \text{ for } y \in U_R,
 \end{equation*} 
 where $C_1$ is a large positive constant independent of $\mu, R.$  
\end{lemma}
\begin{proof}
 Thanks to Lemma \ref{lowbdlem}, we need to prove the lower and upper estimates only when $\mu$ is small. Hence we consider there exists $\widetilde{\mu}>0,$ such that $\mu\in (0,\widetilde{\mu}).$ Let us define 
\begin{equation*}
    \begin{aligned}
        \widetilde{u}_{\mu}(x)= \frac{\umr(Rx)}{\mu R^{s}} \text{ for } x\in U.
    \end{aligned}
\end{equation*}
Then it solves 
\begin{equation}\label{secondlemmatransformationeqn}
    \begin{aligned}
        \fp \widetilde{u}_{\mu} + \mu^{q-p}\fq \widetilde{u}_{\mu}&= 1+\widetilde{C} \mu^{q-p} \text{ in } U,
        \\
       \widetilde{u}_{\mu} &= 0  \text{ in } U^{c}.
    \end{aligned}
\end{equation}
Using $\widetilde{u}_{\mu}$ as a test function in the weak formulation, and the continuous embedding $\wsp \hookrightarrow L^1$ for any $p>1,$ we get $[\;\widetilde{u}_{\mu}\;]_{W_{0}^{s,p}}\leq C,$ independent of $\mu.$
Hence $\widetilde{u}_{\mu}\rightharpoonup v_0$ in $\wsp(U)$ upto a subsequence. Using Theorem \ref{careg} in Appendix, we have that $\|\widetilde{u}_{\mu}\|_{L^{\infty}(U)}$ is uniformly bounded for small values of $\mu.$ Hence, for $s\in (0,\frac{1}{q}),$ we can show that 
    \begin{equation*}
        \begin{aligned}
            \int_{\mathbb{R}^N \times \mathbb{R}^N } \dfrac{(\widetilde{u}_{\mu}(x)-\widetilde{u}_{\mu}(y))^{q-1}(\varphi(x)-\varphi(y))}{|x-y|^{N+sq}}  \,dx\,dy \leq C 
        \end{aligned}
    \end{equation*}
    for each $\varphi \in C_c^\infty(\Om)$ where $C$ is independent of $\mu.$ Clearly as $\mu\rightarrow 0,$ 
     \begin{equation*}
        \begin{aligned}
            \mu^{q-p}\int_{\mathbb{R}^N \times \mathbb{R}^N } \dfrac{(\vM(x)-\vM(y))^{q-1}(\varphi(x)-\varphi(y))}{|x-y|^{N+sq}}  \,dx\,dy \rightarrow 0 \mbox{ for all } \varphi\in C_c^\infty(U).
        \end{aligned}
    \end{equation*}
    Due to the weak-weak continuity property \cite[Lemma 2.2]{CMM18} of fractional $p$-Laplacian and the density argument,  the function $v_0$ solves
    \begin{equation}\label{evo1}
        \begin{aligned}
            \fp v_0 &=1 \ \text{in} \ U,\\
v_0 &= 0 \hspace{0.2cm} \text{in} \ U^{c}.
        \end{aligned}
    \end{equation}
    Thanks to Theorem \ref{careghld}, we can use Ascoli-Arzela theorem to infer that $\widetilde{u}_{\mu}$ has a uniformly convergent subsequence. Uniqueness of the solution of \eqref{evo1} would imply that   $\|\widetilde{u}_{\mu} - v_0\|_{L^\infty(K_1)}\rightarrow0$  as $\mu \to 0$ for any fix compact subset $K_1$ of $U.$ Now, 
   using the strong comparison principle for fractional $p$-Laplacian,  we have  $\underset{K_1}{\inf} \, v_0\geq C_1  >0.$ Since $\widetilde{u}_{\mu}$ uniformly converges to $v_0,$ for sufficiently small $\mu$, we can conclude that 
   \begin{equation*}
    \begin{aligned}
      \widetilde{u}_{\mu}(x)\geq v_0(x)-\varepsilon \geq C_1 -\varepsilon \ \text{for all } \ x\in  K_1.
    \end{aligned}
\end{equation*}
 Now, choosing $\varepsilon$ small enough we obtain, 
\begin{equation}\label{vmueqn}
    \begin{aligned}
        \underset{K_1}{\inf} \, \widetilde{u}_{\mu} \geq C>0  \, \text{  for all } \mu< \widetilde{\mu} \,\, \text{ and } \forall x\in K_1 ,
    \end{aligned}
\end{equation}
where $C=C(N,s,p,\Omega)$ is a constant independent of $\mu.$ Now, let  
\begin{equation*}
    \begin{aligned}
        k=\min_{t=p,q}\left\{\left(\frac{1}{2\max_{U\setminus K_1}( \ft \dis_{U}^{s})}\right)^\frac{1}{t-1}, \frac{C }{\max_{K_1}\dis_{U}^{s}} \right\},
    \end{aligned}
\end{equation*}
where $C>0$ is same as defined in the previous line.
 Now using comparison principle we obtain
 \begin{equation}\label{deonevmu}
     \begin{aligned}
         k \dis_{U}^{s} \leq \widetilde{u}_{\mu} \text{ in } \RR^N.
     \end{aligned}
 \end{equation}  
Now substitute the expression of $\widetilde{u}_{\mu}$ and obtain the required lower bound for all $\mu>0.$
Finally from \eqref{secondlemmatransformationeqn} and Theorem \ref{careg}, the upper estimate follows exactly as in the proof of Lemma \ref{lowbdlem}.
\end{proof}

\section{The Lower Bound}\label{lowbdsec}
In this section, our focus is directed towards the analysis of supersolutions for problems to \eqref{eqn1} on special domains as in \cite{IMS20}. These supersolutions, denoted as $u,$ are assumed to be bounded below by $m \doos.$ Main result of this section is given in Proposition \ref{lowebound}  where we obtain a lower bound for $(\frac{u}{\doos} -m)$ near the boundary in terms of nonlocal excess defined below:
\begin{equation*}
    \operatorname{Ex}(u,m,R,x_0)= \intavg_{\brx} \left|\frac{u(x)}{\doos(x)}-m\right|\thinspace dx.
\end{equation*} 
Let us assume that  $0 \in \partial\Om,$ $R \in (0, \rho/4)$ where $\rho$ is as given in Section \ref{prelim}.  Define 
\begin{equation*}
    A_R:= \bigcup\Bigg\{B_{r}(y)\thinspace : \thinspace y \in \RR^N, \thinspace r \geq \frac{R}{8}, \thinspace B_{r}(y) \thinspace \subset D_R \Bigg\}.
\end{equation*}
Then by equation (3.3) of \cite{IMS20}, we have for some $C>0$
\begin{equation}\label{disone}
    \doo \leq C \dar \text{ in } D_{R/2}.
\end{equation}
Define 
\begin{equation*}
    \begin{aligned}
        A_1:=\{x \,| \, Rx\in A_R\}.
    \end{aligned}
\end{equation*}
Then by definition, we have
\begin{equation}\label{diseqn}
    \begin{aligned}
        R \daone(x) = \dar (Rx) \text{ for } x\in A_1.
    \end{aligned}
\end{equation}
Now we consider a function $u\in\tw(D_R)$ which satisfies the following in the weak sense for some constants $\kp,\kq,\hp,\hq \geq 0$  and $m \geq m_0$ for some fixed constant $m_0>0,$
\begin{eqnarray}\label{subsoleqn}
    \left\{
    \begin{array}{rll}
       \fp u + \fq u \thinspace &\geq \thinspace \sum_{t=p,q}-\kt - m^{t-2}\htt &\text{ in } D_R
        \\
         u \thinspace &\geq \thinspace m \doos   &\text{ in } \RR^N.
    \end{array}
    \right.
\end{eqnarray}
 We now fix the constants $K_t, H_t, m$ for $t=p,q$ through out this section.
\begin{lemma}\label{lwbdlemma}
    Let $u \in\tw(D_R) $ solves \eqref{subsoleqn}. Then there exist constants $\theta_1(N,p,q,s,\Om)\geq 1,$ $C_{3,t}(N,t,s,\Omega)>1$ for $t=p,q$ and $\sigma_1(N,p,q,s,\Om)\in (0,1]$ such that if $\ex\geq m \theta_1$ for all $R\in (0,\rho/4),$ then 
    \begin{equation*}
        \begin{aligned}
           & \inf_{D_{\frac{R}{2}}} \left(\frac{u(x)}{\doos(x)}-m\right) \geq \sigma_1 \ex + \sum_{t=p,q} \left(- C_{3,t} (\kt R^s)^{\frac{1}{t-1}} 
            - C_{3,t} \htt R^s\right).
        \end{aligned}
    \end{equation*}      
\end{lemma}
\begin{proof}
 Let $\Tilde{v}\in \wo(A_R)$ satisfies 
 \begin{eqnarray*}
    \left\{
    \begin{array}{rll}
     \fp \Tilde{v} + \fq \Tilde{v} &= \frac{\lambda^{p-1}+  \lambda^{q-1}}{R^s} \ &\text{in} \ A_R
            \\
           \Tilde{v}&=0 &\text{in} \ A_{R}^{c}.
    \end{array}
    \right.
\end{eqnarray*}
 Using \eqref{disone} and Lemma \ref{lowbdlem},  there exists a positive constant $C>1$ large enough such that $\Tilde{v}(x)\geq \frac{\lambda \doos(x)}{C}$ for $x\in D_{R/2}.$ Henceforth, we will maintain the specific value of $C$ consistently throughout the proof of this lemma. We define $v\in  \wo(A_R)$ satisfying
 \begin{eqnarray*}
     \left\{
     \begin{array}{rll}
           \fp v + \fq v &= \frac{\lp}{R^s} \ &\text{in} \ A_R
            \\
           v&=0 &\text{in} \ A_{R}^{c}.   
     \end{array}\right.
 \end{eqnarray*}
Since $(2C)^{\frac{q-p}{p-1}}>1,$ using comparison principle we get  
\begin{equation}
    \label{disst}
    \begin{aligned}
        v\geq \Tilde{v} \geq \frac{\lambda \doos}{C} \text{ in } D_{R/2}.
    \end{aligned}
\end{equation} We define  $w$ such that 
    \begin{equation*}
        w(x)=
        \left\{
        \begin{aligned}
            &v(x)\ \text{in} \ \bR^c
            \\
            &u(x) \ \text{in} \ \bR.
        \end{aligned}
        \right.
    \end{equation*}
Since $\dis(\bR, D_R)>0$ we can use Proposition \ref{superpostn} of Appendix to infer that 
\begin{equation*}
    \begin{aligned}
         \fp w(x)+ \fq w(x)
           =&\fp v(x) + \fq v(x) \\&+\sum_{t=p,q} 2 \int_{\bR} \frac{(w(x)-u(y))^{t-1}-w^{t-1}(x)}{|x-y|^{N+ts}} \thinspace dy 
    \end{aligned}
\end{equation*}
for $x\in D_R.$  From the calculations of \cite[pp. 20]{IMS20} we get 
\begin{equation*}
    \begin{aligned}
         2 \int_{\bR} \frac{(w(x)-u(y))^{t-1}-w^{t-1}(x)}{|x-y|^{N+ts}} \leq  -\frac{\ex^{t-1}}{CR^s} \text{ for } t=p,q.
    \end{aligned}
\end{equation*}
Thus, we have 
\begin{equation*}
    \begin{aligned}
         \fp w(x)&+ \fq w(x) \leq\frac{\lp}{R^s}-\frac{\ex^{p-1}}{CR^s}-\frac{\ex^{q-1}}{CR^s}.
    \end{aligned}
\end{equation*}
If we choose 
\begin{equation*}
\begin{aligned}
    \lambda=\frac{\ex}{(2C)^{\frac{1}{p-1}}},
\end{aligned}
\end{equation*}
then, we have 
\begin{equation}\label{ew}
    \begin{aligned}
       &\fp w(x)+ \fq w(x) \\& \leq -\frac{\ex^{p-1}}{2CR^s}-\frac{\ex^{q-1}}{2CR^s} \ \text{in} \ D_R .  
    \end{aligned}
\end{equation}
Now we choose 
\begin{equation*}
    \begin{aligned}
       & \theta_1=\frac{1}{\sigma_{1}}=2 C(2 C)^{\frac{1}{p-1}} \geq 1,  \\
& C_{3,t}=\sigma_{1} \max \left\{(4 C)^{\frac{1}{{t-1}}}, 4 C \theta_{1}^{2-t}\right\} \geq 1  \text{ for } t=p,q.
    \end{aligned}
\end{equation*}
Since $\frac{u}{\doos}\geq m$ in $\RR^N,$ the only relevant scenario that requires consideration is the case of   
\begin{equation*}
    \begin{aligned}
        \sigma_1 \ex \geq \sum_{t=p,q} C_{3,t} (\kt R^s)^{\frac{1}{t-1}} +  C_{3,t} \htt R^s.
    \end{aligned}
\end{equation*}
By our choice of $\theta_1$ and $C_{3,t}$ we have for $t=p,q$
 \begin{equation*}
     \ex^{t-1}\left\{
     \begin{aligned}
         &\geq \left(\frac{C_{3,t}}{\sigma_1}\right)^{t-1}\kt R^s \geq 4C\kt R^s
         \\
         & \geq (m\theta_1)^{t-2} \ex \geq (m\theta_1)^{t-2} \frac{C_{3,t}}{\sigma_1} \htt R^s 
        \geq m^{t-2} 4C \htt R^s.
     \end{aligned}\right.
 \end{equation*}
  Summing up we get that
 \begin{equation*}
     \begin{aligned}
         \ex^{t-1} \geq 2CR^{s}(\kt+m^{t-2}\htt).
     \end{aligned}
 \end{equation*}
 Hence using \eqref{subsoleqn} and \eqref{ew}, we have for $x\in D_R$
 \begin{equation*}
     \begin{aligned}
         \fp w(x)+ \fq w(x)  &\leq -\frac{\ex^{p-1}}{2CR^s}-\frac{\ex^{q-1}}{2CR^s} 
        \\  &\leq \sum_{t=p,q}-\kt - m^{t-2}\htt \ 
         \\
         & \leq  \fp u(x)+ \fq u(x).
     \end{aligned}
 \end{equation*}
 Since $w=\chi_{\bR} u$ in $D_{R}^{c},$ we can use comparison principle and \eqref{disst} to conclude that
 \begin{equation*}
     \begin{aligned}
           u(x) \geq \frac{\lambda}{C} \doos(x) = \frac{\ex}{C(2C)^{\frac{1}{p-1}}} \doos(x) \text{ in } D_{R/2}.
     \end{aligned}
 \end{equation*}
 By our assumption $\ex \geq m \theta_1$ and with the choice of $\theta_1$, we obtain 
\begin{equation*}
    \begin{aligned}
        \inf_{D_{\frac{R}{2}}} \left(\frac{u}{\doos}-m\right) \geq \ex\left(\frac{1}{C(2C)^{\frac{1}{p-1}}}-\frac{1}{\theta_1}\right)=\sigma_1 \ex.
    \end{aligned}
\end{equation*}
\end{proof}
In the following lemma, we modify the barrier construction originally designed for the fractional $p$-Laplacian to formulate a barrier for the fractional $(p,q)$-Laplacian.
\begin{lemma} \label{iannlemma}
    For all $\lambda>0,$ we define $w_\lambda(x)= m (1+ \lambda \varphi(\frac{x}{R})) \doos(x)$ for some $\varphi\in C^{\infty}_c(B_1)$ such that $0 \leq \varphi \leq 1,$ $\varphi=1$ in $B_{\frac{1}{2}}.$ Then there exists $C_5(N,p,q,s,\Om)>0$ such that for all $0<\lambda \leq \lambda_0$
    \begin{equation}\label{wlameqn}
        \begin{aligned}
            \fp w_\lambda + \fq w_\lambda \leq C_4(1+\frac{\lambda}{R^s})(m^{p-1}+m^{q-1}) \quad \text{in} \ D_R
        \end{aligned}
    \end{equation}
    
\end{lemma}
\begin{proof}
    Using \cite[Lemma 3.4]{IMS20} for $R$ in the place of $R/2,$  we can guarantee the existence of constants $\lambda_1(N, p, s, \Omega, \varphi)>0, C_5(N, p, s, \Omega, \varphi)>0$ such that for all $|\lambda| \leqslant \lambda_1$
\begin{equation*}
    \begin{aligned}
        \left|(-\Delta)_p^s  w_\lambda\right| \leqslant C_5 m^{p-1}\left(1+\frac{|\lambda|}{R^s}\right) \quad \text { in } D_{R}.
    \end{aligned}
\end{equation*}
Similarly there exist constants like $\lambda_2(N, q, s, \Omega, \varphi)>0,$ and $C_6(N, q, s, \Omega, \varphi)>0$ such that for all $|\lambda| \leqslant \lambda_2$
\begin{equation*}
    \begin{aligned}
        \left|(-\Delta)_q^s  w_\lambda\right| \leqslant C_6 m^{q-1}\left(1+\frac{|\lambda|}{R^s}\right) \quad \text { in } D_{R}.
    \end{aligned}
\end{equation*}
Set, $C_4= \max \{C_5, C_6\}$ and $\lambda_0= \min \{\lambda_1, \lambda_2\}$ to conclude the result.
\end{proof}

For any $\theta_1>1,$ either $\ex \geq m \theta_1$ or $\ex \leq m \theta_1.$ For both the cases, our goal is to establish a lower bound for $\left(\frac{u}{\doos}-m\right).$ If $\ex \geq m \theta_1$ does not hold, then we want to prove the lower bound for any $\theta>1.$ The upcoming lemma is designed to provide this lower bound in such cases.

\begin{lemma}\label{lamlex}
    Let $u \in\tw(D_R) $ solves \eqref{subsoleqn}. Then for all $\theta\geq 1$ there exist $C_{\theta,t}(N,p,q,s,\Omega,\theta)>1,$  for $t=p,q$ and $0<\sigma_\theta(N,p,q,s,\Omega,\theta)\leq 1$ such that if $\ex\leq m \theta$ for all $R\in (0,\rho/4),$ then 
        \begin{equation*}
        \begin{aligned}
           \inf_{D_{\frac{R}{2}}} \left(\frac{u(x)}{\doos(x)}-m\right) \geq \sigma_\theta &\ex 
         +\sum_{t=p,q}\Big(- C_{\theta,t} (m^{t-1}+\kt)^{\frac{1}{t-1}} R^{\frac{s}{t-1}}
           -C_{\theta,t} \htt R^s\Big).
        \end{aligned}
    \end{equation*}
\end{lemma}
\begin{proof}
    For all $\lambda>0,$ we define $w_\lambda(x)= m (1+ \lambda \varphi(\frac{x}{R})) \doos(x)$ for some $\varphi\in C^{\infty}_c(B_1)$ such that $0 \leq \varphi \leq 1,$ $\varphi=1$ in $B_{\frac{1}{2}}.$ 
    For all $x \in \RR^N,$ set 
    \begin{eqnarray*}
       v_\lambda(x)=
    \left\{
    \begin{array}{rll}
        &w_\lambda(x) \ &\text{in} \ \bR^c
            \\
            &u(x) \ &\text{in} \ \bR.
    \end{array}
    \right.
    \end{eqnarray*}
Without any loss of generality we assume $\lambda_0 \leq \min \{1, \frac{(3/2)^s-1}{2}\}$. We have for all $x\in D_R$ 
\begin{equation}\label{wlameqn2}
    \begin{aligned}
         \fp v_\lambda + \fq v_\lambda(x) =& 
        \fp w_\lambda + \fq w_\lambda(x) 
        \\
        &+\sum_{t=p,q} 2 \int_{\bR} \frac{(w_\lambda(x)-u(y))^{t-1}-(w_\lambda(x)-w_\lambda(y))^{t-1}}{|x-y|^{N+ts}} \thinspace dy.
        \\
    \end{aligned}
\end{equation}
Following the calculation of \cite[pp. 31-32]{IMS20}, for $C>0$ large enough,  we get that
\begin{equation}\label{wlameqn3}
    \begin{aligned}
         2 \int_{\bR} \frac{(w_\lambda(x)-u(y))^{t-1}-(w_\lambda(x)-w_\lambda(y))^{t-1}}{|x-y|^{N+ts}} \thinspace dy \leq -\frac{m^{t-2}\ex}{CR^s}.
    \end{aligned}
\end{equation}
Combining \eqref{wlameqn}, \eqref{wlameqn2}, \eqref{wlameqn3}, for all $x  \in D_{R}$  we get
\begin{equation*}
    \begin{aligned}
         \fp v_\lambda + \fq v_\lambda(x) \leq  \sum_{t=p, q} \left[C m^{t-1}+\frac{m^{t-2}}{R^{s}}\left(C \lambda m-\frac{\ex}{C}\right)\right].  
    \end{aligned}
\end{equation*}
Fixing $\theta\geq 1$ we choose 
\begin{equation*}
    \begin{aligned}
         \sigma_{\theta} & =\frac{\lambda_{0}}{2 \theta C^{2}}, \quad  \\
C_{\theta,t} & = \sigma_{\theta} \max \left\{4 C,\left(4 C^{2} \theta^{t-2}\right)^{\frac{1}{t-1}}\right\},
\\
\lambda & =\frac{\sigma_{\theta} \ex}{m},
    \end{aligned}
\end{equation*}
where $\lambda_0$ is defined in Lemma \ref{iannlemma}.
Since $\ex\geq m \theta,$ $\lambda_0\leq 1$ and $\theta\geq 1$, we can conclude that $C \lambda m \leqslant  \frac{\ex}{2 C}$ and $ \lambda  \leqslant \frac{\lambda_{0}}{2 C^{2}} \leqslant \lambda_0.$ Using these estimates we get 
\begin{equation}\label{vleqn1}
    \begin{aligned}
         &\fp v_\lambda + \fq v_\lambda(x) 
        \leq  \sum_{t=p, q} \left[C m^{t-1}-\frac{m^{t-2}}{R^{s}}\frac{\ex}{2C}\right] \ \forall \ x \ \in D_{R}.
    \end{aligned}
\end{equation}
Following the calculations of \cite[pp. 33]{IMS20} we can get
\begin{equation}\label{vleqn2}
    \begin{aligned}
        m^{t-2} \ex \geqslant 2 C R^{s}\left(C{m}^{t-1}+\kt+m^{t-2} \htt \right).
    \end{aligned}
\end{equation}
Combining \eqref{vleqn1}, \eqref{vleqn2} and \eqref{subsoleqn},
\begin{equation*}
    \begin{aligned}
     &\fp v_\lambda + \fq v_\lambda(x) 
        \leq   \sum_{t=p, q} -\kt - m^{t-2}\htt \leq \fp u  + \fq u \text{ in } D_R.   
    \end{aligned}
\end{equation*}
Since 
\begin{eqnarray*}
  v_\lambda(x)=  \left\{
    \begin{array}{rll}
        & m \doos(x), &\text{ in } D_{R}^{c}\cap \bR^c,
        \\
        & u  &\text{ in } \bR,
    \end{array}
    \right.
\end{eqnarray*}
and $m \doos\leq u$ in $\RR^N,$ we can use the comparison principle to conclude $ v_{\lambda} \leqslant u$ in $\RR^N.$ More precisely, due to our careful choice of constants, when $x \in D_{\frac{R}{2}}$ we have 
\begin{equation*}
    \begin{aligned}
        \frac{u}{\doos}-m \geq \frac{w_\lambda}{\doos}-m= m \lambda\geq \sigma_\theta \ex
    \end{aligned}
\end{equation*}
which concludes the proof. 
\end{proof}
Next we localize the global bound from below in \eqref{subsoleqn} and prove the main result
of this section. For some $\Tilde{\kp},\Tilde{\kq} \geq 0$ and $m\geq m_0$ for a fixed constant $m_0>0$, we consider $u \in \two(D_R)$ satisfying
\begin{eqnarray}\label{localeqn1}
    \begin{array}{rll}
         \fp u + \fq u &\geq -\Tilde{\kp}-\Tilde{\kq} \ &\text{in} \ D_R,
        \\
        u& \geq m \doos \ &\text{in} \ D_{2R}.
    \end{array}
\end{eqnarray}

\begin{proposition}\label{lowebound}
Let $u \in \two(D_R)$ solve \eqref{localeqn1} for $m\geq m_0$ for a fixed constant $m_0>0.$ There exist $0<\sigma_2\leq 1,$  $C_6 >1$ depending on $N,p,q,s,\Omega$ and for all $\varepsilon>0, $ a constant $\Tilde{C_\varepsilon}=\Tilde{C_\varepsilon}(N,p,q,s,\Omega,\varepsilon)>0$ such that for all $0<R< \rho/4,$
 we have 
 \begin{equation*}
     \begin{aligned}
         \inf_{D_{R/2}} (\frac{u}{\doos}-m) \geq& \sigma_2 \ex- \varepsilon \Big\|\frac{u}{\doos}-m \Big\|_{L^{\infty}(D_R)}-C_6 \tail_1\Big(\big(-\frac{u}{\doos}+m\big)_{+},2R \Big)R^s
         \\
         &- \sum_{t=p,q} \Tilde{C_\varepsilon} \Bigg[m+\Tilde{K_t}^{\frac{1}{t-1}}+\tail_{t-1}\Big(\big(-\frac{u}{\doos}+m\big)_{+},2R \Big)\Bigg]R^{\frac{s}{t-1}}.
     \end{aligned}
 \end{equation*}
\end{proposition}
\begin{proof}
    Without loss of generality we assume $\frac{u}{\doos}-m\in L^{\infty}(D_R).$ Fix $\varepsilon>0$ and define $v= u\vee m\doos.$ Using Proposition \ref{maxminlem} of Appendix, if we write $\var^{t-1}$ in the place of $\var$ then we get 
    \begin{equation*}
        \begin{aligned}
            \fp v + \fq v \geq \sum_{t=p,q} - K_t- m^{t-2} H_t,
        \end{aligned}
    \end{equation*}
    where we define
    \begin{equation*}
        \begin{aligned}
            K_t&=\Tilde{K_t}+\frac{\var^{t-1}}{R^s}\Big\|\frac{u}{\doos}-m\Big\|_{L^{\infty}(D_R)}^{t-1}+ C_{\var,t} \tail_{t-1} \Big(\big(-\frac{u}{\doos}+m\big)_{+},2R \Big)^{t-1},
            \\
            H_t&=C_{2,t} \tail_1\Big(\big(-\frac{u}{\doos}+m\big)_{+},2R \Big).
        \end{aligned}
    \end{equation*}
    Observe that $v$ satisfies \eqref{subsoleqn}. Thanks to Lemma \ref{lwbdlemma} we can find $0<\sigma_1\leq 1\leq \theta_1,$ and $ C_{3,t}(N,p,q,s)\geq 1$ for $t=p,q$ such that if  $\exv \geq m \theta_1$, then 
    \begin{equation}
        \begin{aligned}
             \inf_{D_{\frac{R}{2}}} \left(\frac{v}{\doos}-m\right) \geq \sigma_1 \exv + \sum_{t=p,q}&\Big(- C_{3,t} (\kt R^s)^{\frac{1}{t-1}} 
        - C_{3,t} \htt R^s\Big).
        \end{aligned}
    \end{equation}
     Otherwise, $\exv\leq m\theta_1$ and we choose $\theta=\theta_1\geq 1$ in Lemma \ref{lamlex}. Now, there exists constants  $C_{\theta_1,t} \geq 1$ for $t=p,q$ such that $0<\sigma_{\theta_1}\leq 1 \leq C_{\theta_1,t}$ and 
      \begin{equation}
        \begin{aligned}
           \inf_{D_{\frac{R}{2}}} \left(\frac{u(x)}{\doos(x)}-m\right) &\geq \sigma_{\theta_1} \exv 
        +\sum_{t=p,q}\Big(- C_{\theta_1,t} (m^{t-1}+\kt)^{\frac{1}{t-1}} R^{\frac{s}{t-1}}
      -C_{\theta,t} \htt R^s\Big).
        \end{aligned}
    \end{equation}
Since $v=u$ in $D_{2R}\supset \bR,$ we get that $\ex=\exv.$
Set $\sigma_2= \min\{\sigma_1, \sigma_{\theta_1}\}<1$ and $C=\underset{t=p,q}{\max} \{C_{3,t},C_{\theta_1,t}\} \geq 1$  to conclude that
\begin{equation*}
    \begin{aligned}
        \inf_{D_{\frac{R}{2}}} \left(\frac{u(x)}{\doos(x)}-m\right) &\geq \sigma_{2} \ex+\sum_{t=p,q}\Big(- C (m^{t-1}+\kt)^{\frac{1}{t-1}} R^{\frac{s}{t-1}}
          -C \htt R^s\Big)
          \\
          &\geq  \sigma_2 \ex- C\var \|(\frac{u}{\doos}-m)\|_{L^{\infty}(D_R)} -C \tail_1((-\frac{u}{\doos}+m)_{+},2R )R^s
         \\
         &\hspace{2.5cm}- \sum_{t=p,q} C  \Big[m+\Tilde{K_t}^{\frac{1}{t-1}}+\Tilde{C_\var}\tail_{t-1}((-\frac{u}{\doos}+m)_{+},2R )\Big]R^{\frac{s}{t-1}}
         .
    \end{aligned}
\end{equation*}
Choosing $\var$ small enough and adjusting the constant we conclude the result.
\end{proof}

\section{The Upper Bound}\label{upbdsec}
Let $u$ denote a subsolution to a problem similar to \eqref{eqn1} defined on a special domain (see \eqref{upbd1}). Moreover, we assume $u$ is locally bounded from above by $M \doos.$ Our objective in this section is to establish an upper bound for $(M-\frac{u}{\doos}).$  Throughout this section, as before we shall assume that  $0 \in \partial\Om,$ $R \in (0, \rho/4).$  
Let us define 
\begin{equation*}
    E_R:= \bigcup\Bigg\{B_{r}(y)\thinspace : \thinspace y \in \Om, \thinspace r \geq \frac{R}{8}, \thinspace B_{r}(y) \thinspace \subset D_{4R}\setminus D_{3R/4} \Bigg\}.
\end{equation*}
By (4.2) of \cite{IMS20}, we have for some $C>0$
\begin{equation}\label{dissone}
    \doo \leq C \der \text{ in } D_{3R}\setminus D_R.
\end{equation}
Define 
\begin{equation*}
    \begin{aligned}
        E_1:=\{x| Rx\in E_R\}.
    \end{aligned}
\end{equation*}
Then by definition, we have
\begin{equation}\label{disseqn}
    \begin{aligned}
        R \deone(x) = \der (Rx) \text{ for } x\in E_1.
    \end{aligned}
\end{equation}
\begin{lemma} \label{lemdoos}
    Given an $\kappa>0,$ there exists a $C_1>1$ large enough such that if $\va\in \wo(E_R)$ satisfy
    \begin{eqnarray*}
        \begin{array}{rll}
       \fp\va+\fq\va&=\frac{C_1}{R^s} &\text{ in } E_R,
            \\
            \va&=0 &\text{ in } E^{c}_R,
        \end{array}
    \end{eqnarray*}

 then  
    \begin{equation*}
        \begin{aligned}
            \va(x)\geq \kappa \doos \text{ for all } x\in D_{3R}\setminus D_R.
        \end{aligned}
    \end{equation*}
\end{lemma}
\begin{proof}
    Set  $\va_1(x)=\va(Rx)$ for $x\in E_1.$ Then $\va_1$ would satisfy
    \begin{equation*}
    \begin{aligned}
        R^{(q-p)s} \fp \va_1+ \fq \va_1&=C_1 R^{(q-1)s} \ \text{in} \  E_1,\\
\va_1(x) &= 0 \hspace{1.5cm} \text{in} \ E_{1}^{c}.
    \end{aligned}
    \end{equation*}
   where $C_1>1$ is a constant that would be chosen later. Defining  $w(x)= k R^s \deones(x)$ for $x \in E_1,$ we get 
    \begin{equation*}
        \begin{aligned}
     R^{(q-p)s} \fp w+ \fq w&= R^{(q-1)s}\left(k^{p-1} \fp \deones+ k^{q-1} \fq \deones\right).         
        \end{aligned}
    \end{equation*}
  Now we have $ \underset{E_1\setminus E_K}{\max} \{|\fp \deones|, |\fq \deones| \} \leq M_1,$
    where $M_1$ is a positive constant and $E_K$ is a compact subset of $E_1.$ 

      If we set $v(x)=\frac{\va_1(x)}{R^s}$ for $x \in E_1$ then $v$ satisfies  
      \begin{eqnarray*}
          \begin{array}{rll}
                     \fp v + \fq v &=C_1 \ &\text{in} \ E_1,
         \\
         v&=0  \ &\text{in} \ E_{1}^{c}.
          \end{array}
      \end{eqnarray*}

     Thanks to  \eqref{Ceqn}, we can infer that 
    \begin{equation*}
        \begin{aligned}
       \inf_{E_K}    \frac{\va_1(x)}{R^s} = \inf_{E_K} v \geq C_2 (C_1)^{\frac{1}{q-1}} >0 \text{ for all } &x\in E_K,
       \text{ where } C_2 \text{ does not depend on } C_1.
        \end{aligned}
    \end{equation*}
    Choose  
    \begin{equation}
        \begin{aligned}
            k = \min  \left\{\left(\frac{C_1}{2M_1}\right)^{1/p-1},\left(\frac{C_1}{2M_1}\right)^{1/q-1},\frac{C_2(C_1)^{\frac{1}{q-1}}}{\underset{E_1}{\max}\deones} \right\},
        \end{aligned}
    \end{equation}
    to use the comparison principle on 
    \begin{eqnarray*}
        \left\{
        \begin{array}{rll}
          &R^{(q-p)s} \fp \va_1+ \fq \va_1=C_1 R^{(q-1)s} \geq   R^{(q-p)s} \fp w+ \fq w \ &\text{in} \ E_1\setminus E_K,
               \\
           &\va_1 \geq \underset{E_K}{\inf} \va_1 \geq  C_2 (C_1)^{\frac{1}{q-1}} R^s \geq  k R^s \deones=w(x) \hspace{1.8cm} &\text{in} \ E_K.
        \end{array}\right.
    \end{eqnarray*}

    Combining with \eqref{dissone} we get that 
    \begin{equation*}
    \begin{aligned}
        \va(y) \geq k \ders(y) \geq k C_5 \doos(y) \text{ for } y\in D_{3R}\setminus D_R.
    \end{aligned}
\end{equation*}
Now given $\kappa, M_1,C_2, C_5$ fixed, we can choose $C_1$ large enough such that 
 \begin{equation*}
        \begin{aligned}
          \frac{\kappa}{C_5}\leq  k = \min \left\{\left(\frac{C_1}{2M_1}\right)^{1/p-1},\left(\frac{C_1}{2M_1}\right)^{1/q-1},\frac{C_2(C_1)^{\frac{1}{q-1}}}{\underset{E_1}{\max}\deones} \right\},
        \end{aligned}
    \end{equation*}
    and hence the result.
\end{proof}

We now construct the barrier as a solution of a double obstacle problem.
\begin{lemma}\label{barlem1}
  Let  $\Bar{x}\in D_{R/2},$   and $R\in(0,\rho/4).$ Fix an $ \widetilde{M}_0>0$ and $M > \widetilde{M}_0,$  there exists a function $v\in \wo(\Om)\cap C(\RR^N)$ and a positive constant $C(N,p,q,s,\Omega)>1$ such that the following conditions are satisfied: 
    \begin{itemize}
        \item[(i)] $|\fp v + \fq v|\leq \frac{C(M^{p-1}+M^{q-1})}{R^s}$ in $D_{2R}$;
         \item[(ii)] $v(\Bar{x})=0;$ 
          \item[(iii)] $v\geq M\doos$ in $D^{c}_R;$
           \item[(iv)] $|v|\leq C R^s $ in $D_{2R}.$
    \end{itemize}
\end{lemma}
\begin{proof}
    First, we construct the lower obstacle.  Define $\va \in \wo(E_R)$ such that 
    \begin{eqnarray*}
        \begin{array}{rll}
           \fp \va+ \fq \va&= \frac{C_1}{R^s} \ &\text{in} \ E_R,\\
\va &= 0 \ &\text{in} \ E_{R}^{c} ,
        \end{array}
    \end{eqnarray*}
where $C_1>0.$ Thanks to Lemma \ref{lemdoos}, for the given $M>0$  we  choose $C_1$ large enough to get 
      \begin{equation} \label{lowobsp3}
        \begin{aligned}
            \va(x) \geq M \doos(x) \text{ for } x\in D_{3R}\setminus D_R.
        \end{aligned}
    \end{equation}
   Since  $\va=0$ in $E^{c}_R,$ and by Lemma \ref{lowbdlem}, we infer that
    \begin{equation}\label{lowobsp2}
        \begin{aligned}
            \va\leq C_2 R^s \text{ in } \RR^N.
        \end{aligned}
    \end{equation}
Given an $M > \widetilde{M}_0>0,$ we can choose $C_2$ large enough such that $C_1 \leq C_2(M^{p-1}+M^{q-1}).$  Using Lemma \ref{extlem}, we have 
    \begin{equation}\label{lowobsp1}
        \begin{aligned}
            \fp \va + \fq \va \leq \frac{C_1}{R^s} \leq \frac{C_2(M^{p-1}+M^{q-1})}{R^s} \ \text{in} \ \RR^N.
        \end{aligned}
    \end{equation}
     Now, fix these constants $C_1$ and $C_2$ in this proof. The function $\va$ constructed here will serve as the lower obstacle in the latter part of the proof. Next, we shall discuss the upper obstacle.
    \\
    Let $\Psi \in \wo(B_{R/8})$ satisfies
    \begin{eqnarray*}
        \begin{array}{rll}
        \fp \Psi + \fq \Psi&= \frac{\nn}{R^s} \ &\text{in} \  B_{R/8},\\
\Psi &= 0 \hspace{2cm} &\text{in} \ B_{R/8}^{c}. 
        \end{array}
    \end{eqnarray*}
Since $\fp$ and $\fq$ both are rotation invariant by definition,  $\Psi$ is radially decreasing. We define
\begin{equation*}
\begin{aligned}
        \psi(x)= \max_{\RR^N} \Psi- \Psi(x-\x) \text{ for }x\in\RR^N. 
\end{aligned}
\end{equation*}
Clearly $\psi \in \tw(\Om),$ $\psi\geq 0.$  Since $\Psi$ is radially decreasing and maximum attained at $0$ thus $\psi(\x)=0.$ Moreover we want to prove 
\begin{equation*}
    \begin{aligned}
            \va \leq \psi \text{ in } \RR^N \text{ for } \nu \text{ large enough.}
    \end{aligned}
\end{equation*}
For $\nu$ large enough, we can employ Lemma \ref{lowbdlem} to get  
\begin{equation*}
    \begin{aligned}
            \Psi(x) \geq k \nu R^s \dis_{B_{\frac{1}{8}}}^{s}(x) \text{ for } x\in \RR^N,
    \end{aligned}
\end{equation*}
for some $k>0.$ Hence we have that 
\begin{equation}\label{max}
    \begin{aligned}
         \max_{\RR^N} \Psi \geq k \nu R^s  \max_{\RR^N} \dis_{B_{\frac{1}{8}}}^{s}(x) \geq \frac{\nu R^s}{C_3} \text{ for } C_3 \text{ large enough.}
    \end{aligned}
\end{equation}
Now to compare $\psi$ and $\varphi,$ we first fix any $x \in \RR^N.$ If $x \in D_{3R/4}$ then , by construction,
\begin{equation*}
    \begin{aligned}
        \va(x)=0\leq \psi(x) \text{ in } D_{3R/4}.
    \end{aligned}
\end{equation*}
If $x \in D^{c}_{3R/4},$  then $|x-\x|>\frac{R}{8}.$  Thus $\Psi(x-\x)=0$ for $x\in D^{c}_{3R/4}.$ We combine \eqref{lowobsp2},\eqref{max} and choose $\nu$ large enough such that $\nu \geq C_3 C_2$ to get 
\begin{equation*}
    \begin{aligned}
        \va(x)\leq C_2 R^s \leq \psi(x) \text{ in }  D^{c}_{3R/4}.
    \end{aligned}
\end{equation*}
We fix this particular choice of $\nu$ for the remaining part of the proof. Again using lemma \ref{lowbdlem}, we have 
\begin{equation}\label{upobsp2}
    \begin{aligned}
        \psi(x)\leq \max_{\RR^N} \Psi \leq C_4 R^s.
    \end{aligned}
\end{equation}
Also we note that using Lemma \ref{extlem} we have 
\begin{equation}\label{upobsp1}
    \begin{aligned}
        \fp \psi +\fq \psi \geq \frac{-C_5}{R^s} \text{ in } \RR^N.
    \end{aligned}
\end{equation}
Ultimately we shall now construct the barrier. 
We know that there exists a unique $\Phi\in \wo(\Omega)$ which minimizes the quantity $\frac{1}{p}[u]^{p}_{s,p}+\frac{1}{q}[u]^{q}_{s,q}$ in the set $\left\{ u\in \wo(\Omega): \va \leq u \leq \psi \text{ in } \RR^N\right\}.$ And, $\Phi$ satisfies, 
    
    \begin{equation*}
        \begin{aligned}
            0 \wedge \left(\fp\psi +\fq \psi\right) \leq \fp \Phi +\fq \Phi \leq 0 \vee \left(\fp\va +\fq \va\right) \text{ in } \Om.
        \end{aligned}
    \end{equation*}
    See Lemma \ref{LSlem} in the Appendix for proof. Thus, for $M>\widetilde{M}_0>0$ and  $C'>0$ large enough, 
    \begin{itemize}
        \item[(a)] $|\fp\Phi +\fq \Phi| \leq \frac{C'(M^{p-1}+M^{q-1})}{R^s} $ in $D_{2R}$ due to  \eqref{lowobsp1} and \eqref{upobsp1};
         \item[(b)] $\Phi(\x)=0$ since $\va(\x)=0=\psi(\x)$; 
          \item[(c)] $   0\leq \Phi \leq \psi \leq C' R^s \text{ in } \RR^N$ thanks to \eqref{upobsp2};
           \item[(d)] $\Phi \geq \va \geq M\doos \text{ in } D_{3R}\setminus D_R$ using \eqref{lowobsp3}.
    \end{itemize}
   We require the property (d) to be satisfied in $D_{R}^{c}$ itself for which we define
   \begin{eqnarray*}
       v(x)=\left\{
       \begin{array}{rll}
        &\Phi(x) \thickspace &\text{ if } x\in D_{3R},
        \\
        &\Phi(x) \vee M\doos(x) \thickspace &\text{ if } x\in D^{c}_{3R}.
       \end{array}
       \right.
   \end{eqnarray*}
Clearly $v\in \wo(\Om)$ satisfies (ii),(iii), and (iv) by construction. It remains to prove the property (i) for the function $v.$ Clearly,  
\begin{eqnarray*}
    \begin{array}{rll}
  \fp v(x)+ \fq v(x)
  =&\fp \Phi(x) + \fq \Phi(x)
  \\
  &+\sum_{t=p,q} 2 \int_{D^{c}_{3R}\cap \{\Phi<M\doos\}} \frac{(\Phi(x)-M\doos)^{t-1}-(\Phi(x)-\Phi(y))^{t-1}}{|x-y|^{N+ts}}.
    \end{array}
\end{eqnarray*}
Due to the monotonicity of the function $r \mapsto r^{t-1}$, the integrand  in the above expression is negative. Consequently, in accordance with property (a) of $\Phi$, this implies
 \begin{equation*}
        \begin{aligned}
            \fp v(x) + \fq v(x)\leq \frac{C'(M^{p-1}+M^{q-1})}{R^s} \text{ in } D_{2R}.
        \end{aligned}
    \end{equation*}
    and gives the upper bound for $v$ as required in $(i).$
Next, we modify the calculation outlined in \cite[pp. 39]{IMS20}, to get that 
   \begin{eqnarray*}
       \begin{array}{rll}
            \sum_{t=p,q} 2 \int_{D^{c}_{3R}\cap \{\Phi<M\doos\}} \frac{(\Phi(x)-M\doos)^{t-1}-(\Phi(x)-\Phi(y))^{t-1}}{|x-y|^{N+ts}} \thinspace dy \geq \int_{D^{c}_{3R}} \frac{-C^{\prime}}{|y|^{N+s}} \thickspace dy,
       \end{array}
   \end{eqnarray*}
   for all  $x\in D_{2R}$ and $y\in D^{c}_{3R}.$ We combine the last inequality with the property (a) of $\Phi$ and obtain the required lower bound for $\fp v(x) + \fq v(x)$ in $D_{2 R}.$ This completes the proof. 
\end{proof}
Let  $u \in \tw(D_R)$ satisfies 
\begin{eqnarray}\label{upbd1}
    \left\{
    \begin{array}{rll}
    \fp u + \fq u &\leq \sum_{t=p,q} \kt + M^{t-2}\htt \ &\text{in} \ D_R,
        \\
        u &\leq M \doos \ &\text{in } \RR^N,
    \end{array}
    \right.
\end{eqnarray}
for some $\kp,\kq,\hp,\hq, \geq 0$ and $M>\widetilde{M}_0>0.$
\begin{lemma}\label{help}
  Let $u \in \tw(D_R)$ satisfy \eqref{upbd1}. Then there exists $C_{4,t}(N,t,s,\Omega)>1$ for $t=p,q$ such that 
    \begin{equation*}
        \begin{aligned}
            \Ex\geq \sum_{t=p,q} C_{4,t} (M+(\kt R^s)^{1/{t-1}} +\htt R^s) ,
        \end{aligned}
    \end{equation*}
    implies that
    \begin{equation*}
        \begin{aligned}
            \sup_{D_{R/2}} u \leq 0.
        \end{aligned}
    \end{equation*}
\end{lemma}
\begin{proof}
   Fix $\x\in D_{R/2}.$ Let us define 
    \begin{equation*}
        w(x)=\left\{
        \begin{aligned}
            & v(x) \text{ in } {\bR}^{c},
            \\
            & u(x) \text{ in } \bR,
        \end{aligned}\right.
    \end{equation*}
     where $v$ is defined in Lemma \ref{barlem1}. Then for all $x \in D_R$
     \begin{equation*}
         \begin{aligned}
              \fp w(x) + \fq w(x)=& \fp  v(x) + \fq  v(x) 
            \\ & + \sum_{t=p,q} 2 \int_{\bR} \frac{(v(x)-u(y))^{t-1}-( v(x)- v(y))^{t-1}}{|x-y|^{N+ts}}.
         \end{aligned}
     \end{equation*}
  Using Lemma \ref{barlem1}, property $(i)$ of $v$  we get 
      \begin{equation*}
         \begin{aligned}
              \fp w(x) &+ \fq w(x)=  \sum_{t=p,q} \frac{-CM^{t-1}}{R^s} +2 \int_{\bR} \frac{(v(x)-u(y))^{t-1}-( v(x)- v(y))^{t-1}}{|x-y|^{N+ts}}.
         \end{aligned}
     \end{equation*}
     Since $v$ satisfies the property $(iii)$ of Lemma \ref{barlem1}, we can  use  a similar calculation of \cite[pp.40]{IMS20} to conclude that 
     \begin{equation*}
         \begin{aligned}
            2 \int_{\bR} \frac{(v(x)-u(y))^{t-1}-( v(x)- v(y))^{t-1}}{|x-y|^{N+ts}} \geq \frac{1}{C} \frac{\Ex^{t-1}}{R^s},
         \end{aligned}
     \end{equation*}
     for $C$ large enough. Next choose $ C_{4,t} \geq (3C^2)^{1/{t-1}} \geq (3C)^{1/{t-1}}$ to get that 
        \begin{equation*}
       \Ex^{t-1}\geq\left\{ 
       \begin{aligned}
           &(C_{4,t} M)^{t-1} \geq 3C^2 M^{t-1}
           \\
           & (C_{4,t})^{t-1} \kt R^s \geq 3 C \kt R^s
           \\
           & (C_{4,t} M)^{t-2} \Ex \geq 3C M^{t-2} \htt R^s
       \end{aligned}
       \right.
   \end{equation*}
   and consequently 
     \begin{equation*}
       \begin{aligned}
           \Ex^{t-1} \geq C^2 M^{t-1} + C \kt R^s + C M^{t-2} \htt R^s.
       \end{aligned}
   \end{equation*}
   So we have for all $x\in D_{R}$
   \begin{equation*}
       \begin{aligned}
         \fp w(x) + \fq w(x) \geq  \sum_{t=p,q}   \kt +  M^{t-2} \htt \geq  \fp u(x) + \fq u(x)  .
       \end{aligned}
   \end{equation*}
 Using Lemma \ref{barlem1}, Property $(iii)$ we have $w(x)=v(x)\geq M \doos(x)\geq u(x)$ inside $D^{c}_{R}\setminus \bR.$ By definition, $u=w$ in $\bR.$ Using comparison principle we get that $  u(\x)\leq w(\x)=0 \text{ for any arbitrary } \x\in D_{R/2}.$  Consequently, $\sup_{D_{R/2}} u \leq 0$, as the choice of $\x\in D_{R/2}$ is arbitrary.
\end{proof}
    Now given any $\Theta_1>1,$  either $\Ex \geq M \Theta_1$ or $\Ex \leq M \Theta_1.$ Similar to Section \ref{lowbdsec}, for both cases we want to prove upper bounds for the subsolutions.
\begin{lemma} \label{usmlex}
Let $u \in \tw(D_R)$ satisfy \eqref{upbd1}. Then there exist the constants $\Theta_1(N,p,q,s,\Om)\geq 1,$ $ 0<\sigma_3(N,p,q,s,\Om)\leq 1,$ $C_{5,t}(N,t,s,\Om)>1,$ for $t=p,q$ such that $\Ex\geq M\Theta_1$  for all $R\in (0,\rho/4)$ then
    \begin{equation*}
        \begin{aligned}
        \inf_{D_{\frac{R}{4}}} \left(M-\frac{u(x)}{\doos(x)}\right) \geq \sigma_3 \Ex + \sum_{t=p,q}\left(- C_{5,t} (\kt R^s)^{\frac{1}{t-1}} -C_{5,t} \htt R^s\right).
        \end{aligned}
    \end{equation*}
\end{lemma}
\begin{proof}
   We set 
   \begin{equation*}
       \begin{aligned}
           H_R:= \bigcup \left\{B_{r}(y): y \in D_{3R/8}, r \geq \frac{R}{16}, B_{r}(y) \subset D_{3R/8})\right\}.
       \end{aligned}
   \end{equation*}
   Then by \cite[Equation (4.16)]{IMS20}, we have for some $C>0$ 
\begin{equation}\label{dishr}
    \doo \leq C \dis_{H_{R}} \text{ in } D_{R/4}.
\end{equation}
Let $\va\in \wo(H_R)$ satisfies 
\begin{eqnarray*}
    \begin{array}{rll}
             \fp \va + \fq \va &= \frac{\lo}{R^s} &\text{ in } H_R,
            \\
            \va&=0 &\text{ in } H^{c}_R.
    \end{array}
\end{eqnarray*}
    Then similar to Lemma \ref{lowbdlem} we can conclude that 
    \begin{equation}\label{valem}
        \begin{aligned}
           \frac{\lambda}{C} \dis_{H_{R}}^{s} \leq  \va \leq C \lambda R^s \text{ in } H_R \text{ for } C \text{ independent of } \lambda,R.
        \end{aligned}
    \end{equation}
    Define 
    \begin{eqnarray*}
         v(x)=\left\{
         \begin{array}{rll}
           &-\va &\text{ in } D_{R/2},
            \\
            & M \doos &\text{ in } D^{c}_{R/2}.
         \end{array}
         \right.
    \end{eqnarray*}
    Clearly  $\va\in \tw(H_R)$ and $\dis(D^{c}_{R/2},H_R)>0.$ Then, 
    \begin{equation*}
        \begin{aligned}
             \fp v + \fq v &=  \fp(- \va) + \fq (-\va)  
         \\&+\sum_{t=p,q} 2 \int_{D^{c}_{R/2}} \frac{(-\va(x)- M\doos(y))^{t-1}- (-\va(x))^{t-1}}{|x-y|^{N+ts}} \thinspace dy
             \\
             &
             \geq -\frac{\lo}{R^s}-\sum_{t=p,q} C \int_{D^{c}_{R/2}} \frac{ (\va(x))^{t-1} +(M\doos(y))^{t-1} }{|x-y|^{N+ts}} \thinspace dy,
        \end{aligned}
    \end{equation*}
     for all $x \in H_R \subset  D_{R/2}$ and for some $C$ depends only on $t=p,q.$ Now for $x\in H_R$ and $y \in B^{c}_{R/2}$ we have that $ C|y-x|>|y|.$ By definition $ \doo(y)\leq |y|.$ Hence using \eqref{valem} we get that 
      \begin{equation*}
        \begin{aligned}
           \int_{D^{c}_{R/2}} \frac{ (\va(x))^{t-1} +(M\doos(y))^{t-1} }{|x-y|^{N+ts}} \thinspace dy &\leq  \int_{B^{c}_{R/2}} \frac{ (C \lambda R^s)^{t-1} +(M|y|^{s})^{t-1} }{C |y|^{N+ts}} \thinspace dy
        \\
        & \leq C (\lambda^{t-1}+M^{t-1}){ \int_{B^{c}_{R/2}} \frac{ R^{s(t-1)} +|y|^{s(t-1)} }{ |y|^{N+ts}}} \thinspace dy
        \\
        & \leq C\frac{(\lambda^{t-1}+M^{t-1})}{R^s},
        \end{aligned}
    \end{equation*}
    for $t=p,q$ and $C>1$ large enough. Exploiting a slight abuse of notation for a positive constant $C$, we obtain the following for $x \in H_R$:
    \begin{equation}\label{veqn1}
        \begin{aligned}
          \fp v + \fq v \geq \sum_{t=p,q} -C\frac{(\lambda^{t-1}+M^{t-1})}{R^s}.
        \end{aligned}
    \end{equation}
     Set 
    \begin{equation*}
        w(x)=\left\{
        \begin{aligned}
            &v(x) \text{ in } \bR^c
            \\
            &u(x) \text{ in } \bR.
        \end{aligned}
        \right.
    \end{equation*}
  Then using Proposition \ref{superpostn} of appendix, we have $w \in \tw(H_R)$ and 
  \begin{equation}\label{veqn2}
      \begin{aligned}
            \fp w +\fq w= &\fp v +\fq v 
        \\ & + \sum_{t=p,q} 2\int_{\bR}\frac{(v(x)-u(y))^{t-1}-(v(x)-M\doos(y))^{t-1}}{|x-y|^{N+ts}} \thinspace dy.
      \end{aligned}
  \end{equation}
  From the calculations of \cite[pp. 43]{IMS20} we get that for $t=p,q$
  \begin{equation}\label{veqn3}
      \begin{aligned}
          2\int_{\bR}\frac{(v(x)-u(y))^{t-1}-(v(x)-M\doos(y))^{t-1}}{|x-y|^{N+ts}} \thinspace dy \geq \frac{\ex^{t-1}}{CR^s}.
      \end{aligned}
  \end{equation}
  for $C>1$ large enough. Combining \eqref{veqn1}, \eqref{veqn2}, \eqref{veqn3} we get that 
  \begin{equation}\label{weqn}
      \begin{aligned}
           \fp w +\fq w  \geq \sum_{t=p,q}\Bigg[ -C\frac{(\lambda^{t-1}+M^{t-1})}{R^s} + \frac{\Ex^{t-1}}{CR^s}\Bigg].
      \end{aligned}
  \end{equation}
  for $C>1$ large enough. Now we fix that chosen $C>1$ which depends on $N,p,q,s,\Om.$ For $C_{4,t}$ is as given in Lemma \ref{help}, we can fix the constants as 
  \begin{equation*}
      \begin{aligned}
           \lambda&=\frac{\Ex}{(4C^2)^{\frac{1}{p-1}}};
            \\
           \Theta_1&= \max_{t=p,q} \left\{2C_{4,t},  \left[\frac{1}{2C^2}-\frac{1}{(4C^2)^{\frac{q-1}{p-1}}}\right]^{\frac{-1}{q-1}},(4C^2)^\frac{1}{p-1}\right\};
           \\
           \sigma_3& =\frac{1}{C(4C^2)^\frac{1}{p-1}};
           \\
           C_{5,t} &= \sigma_3 \max_{t=p,q}\left\{2C_{4,t},(4C)^{\frac{1}{t-1}},\frac{4C}{(\Theta_1)^{t-2}}\right\}.
      \end{aligned}
  \end{equation*}
   By the choice of  $\Theta_1$ we have
    \begin{equation*}
 \frac{1}{2C^2}   \leq \left\{
        \begin{aligned}
 &\left[ \frac{1}{C^2}-\frac{1}{4C^2}-\frac{1}{(\Theta_1)^{p-1}} \right] 
 \\
 & \left[ \frac{1}{C^2}-\frac{1}{(4C^2)^\frac{q-1}{p-1}}-\frac{1}{(\Theta_1)^{q-1}} \right].
        \end{aligned}\right.
    \end{equation*}
    By the choice of $\lambda$ and $\Theta_1$ we have
        \begin{equation}\label{weqn1}
        \begin{aligned}
            \fp w +\fq w &\geq \frac{C}{R^s} \Bigg[\frac{\Ex^{p-1}}{C^2}+\frac{\Ex^{q-1}}{C^2 }
            \\
            &\hspace{2cm}-\frac{\Ex^{p-1}}{4C^2}
            -\frac{\Ex^{q-1}}{(4C^2)^\frac{q-1}{p-1}}
            \\
            &\hspace{4.5cm}-\left(\frac{\Ex}{\Theta_1}\right)^{p-1} -\left(\frac{\Ex}{\Theta_1}\right)^{q-1}
            \Bigg]
            \\
            &
            \geq \frac{\Ex^{p-1}C}{R^s} \Bigg[ \frac{1}{C^2}-\frac{1}{4C^2}-\frac{1}{(\Theta_1)^{p-1}} \Bigg] 
            \\
            &
            \hspace{4cm}+ \frac{\Ex^{q-1}C}{R^s} \Bigg[ \frac{1}{C^2}-\frac{1}{(4C^2)^\frac{q-1}{p-1}}-\frac{1}{(\Theta_1)^{q-1}} \Bigg]
            \\
            & \geq 
            \frac{\Ex^{p-1}}{2 C R^s}  + \frac{\Ex^{q-1}}{2 C R^s}.
        \end{aligned}
    \end{equation}
    Since $\frac{u}{\doos}\leq M$  in $\RR^N,$ only nontrivial case to be considered is  
    \begin{equation*}
        \begin{aligned}
             \sigma_3 \Ex - \sum_{t=p,q} C_{5,t} (\kt R^s)^{\frac{1}{t-1}} -C_{5,t} \htt R^s \geq 0.
        \end{aligned}
    \end{equation*}
    So we have 
    \begin{equation*}
        \Ex^{t-1}\geq \left\{
        \begin{aligned}
           & \left(\frac{C_{5,t}}{\sigma_3}\right)^{t-1} \kt R^s \geq 4C \kt R^s
            \\
            & (M\Theta_1)^{t-2} \frac{C_{5,t}}{\sigma_3}\htt R^s\geq 4CM^{t-2} \htt R^s,
        \end{aligned}\right.
    \end{equation*}
    thus for $t=p,q$ 
    \begin{equation*}
        \begin{aligned}
            \frac{\Ex^{t-1}}{2CR^s}\geq \kt + \htt R^s.
        \end{aligned}
    \end{equation*}
    Substituting this into equation \eqref{weqn1}, we obtain:
     \begin{equation*}
        \begin{aligned}
              \fp w +\fq w& \geq \sum_{t=p,q} \kt + \htt R^s \geq    \fp u +\fq u \text{ in } H_R.
        \end{aligned}
    \end{equation*}
   Now consider $x\in D_{R/2}\cap H^{c}_R,$ then we have
        \begin{equation}
            \Ex\geq\left\{
            \begin{aligned}
                &\frac{C_{5,t}}{\sigma_3} (\kt R^s)^\frac{1}{t-1}+ \frac{C_{5,t}}{\sigma_3}\htt R^s \geq 2 C_{4,t}(\kt R^s)^\frac{1}{t-1}
        + 2 C_{4,t} \htt R^s
                \\
                & M\Theta_1 \geq 2 C_{4,t} M.
            \end{aligned}\right.
        \end{equation} 
       In other words, we get 
        \begin{equation*}
            \begin{aligned}
               \Ex \geq \sum_{t=p,q} C_{4,t}\left(M+(\kt R^s)^\frac{1}{t-1}+\htt R^s\right). 
            \end{aligned}
        \end{equation*} 
        Using Lemma \ref{help}, we get that $\underset{D_{R/2}}{\sup \thinspace}u \leq 0.$ Now considering $x \in H_R^c$, we can delineate our analysis into three distinct scenarios:
    \begin{itemize}
        \item[(a)] if $x\in \bR$ then $w(x)=u(x);$
        \item[(b)] if $x\in D^{c}_{R/2}\cap \bR^c$ then $w(x)=M\doos(x)\geq u(x);$
        \item[(c)] if $x\in D_{R/2}\cap H^{c}_R,$ then we have $w(x)=0 \geq \underset{D_{R/2}}{\sup \thinspace}u \geq u(x).$
        \end{itemize}
         Therefore using the comparison principle on $u$ and $w$ we have $u\leq w$ in $\RR^N.$ Using \eqref{valem} and \eqref{dishr}, we have 
    \begin{equation*}
        \begin{aligned}
            u(x)\leq w(x)=v(x)=-\va(x) &\leq -\frac{\lambda}{C} \doos(x)
            \leq-\frac{\Ex}{C(4C^2)^\frac{1}{p-1}} \doos(x)
            \\
            &=-\sigma_3\Ex\doos(x) \text{ in } D_{R/4}.
        \end{aligned}
    \end{equation*}
    Then we can conclude that
    \begin{equation*}
        \begin{aligned}
            \inf_{D_{R/4}}\left(M-\frac{u}{\doos}\right)\geq -\sup_{D_{R/4}}\frac{u}{\doos} \geq \sigma_3 \Ex.
        \end{aligned}
    \end{equation*}
\end{proof}
Analogous to Lemma \ref{iannlemma}, we can deduce the following lemma, and therefore choose to skip its proof.
\begin{lemma} \label{distlwbd}
     We define $w_\lambda(x)= M (1- \lambda \varphi(\frac{x}{R})) \doos(x)$ for all $\lambda>0,$ and for some $\varphi\in C^{\infty}_c(B_1)$ such that $0 \leq \varphi \leq 1$ in $B_1$ and $\varphi=1$ in $B_{\frac{1}{2}}.$ Then there exists $C_7(N,p,q,s,\Omega)$ such that for all $0<\lambda \leq \lambda_0$
    \begin{equation*}
        \begin{aligned}
            \fp w_\lambda + \fq w_\lambda \geq -C_7(1-\frac{\lambda}{R^s})(M^{p-1}+M^{q-1}) \quad \text{in} \ D_R.
        \end{aligned}
    \end{equation*}
   
\end{lemma}
The proof of the forthcoming lemma closely resembles that of Lemma \ref{lamlex}. With Lemma \ref{distlwbd} now established, guided by the approach detailed in \cite[Lemma 4.4]{IMS20}, one can proceed to demonstrate the subsequent lemma. Therefore, we opt to exclude the detailed proof here.
\begin{lemma}\label{ubgex}
   Let $u \in \tw(D_R)$ satisfy \eqref{upbd1} and $R\in (0,\rho/4)$. Then for all $\Theta \geq 1$ there exist constants $0<\sigma_\Theta(N,p,q,s,\Omega,\Theta)\leq 1,$ $C_{\Theta,t}(N,p,q,s,\Omega,\Theta)>1,$ for $t=p,q$  such that $\Ex\leq M \Theta $  then
    \begin{equation*}
        \begin{aligned}
        \inf_{D_{\frac{R}{2}}} \left(M-\frac{u(x)}{\doos(x)}\right) \geq \sigma_3 \Ex + \sum_{t=p,q}\left(- C_{5,t} (\kt R^s)^{\frac{1}{t-1}} -C_{5,t} \htt R^s\right).
        \end{aligned}
    \end{equation*}  
\end{lemma}
Finally, we establish the counterpart to Proposition \ref{lowebound}. We deal with $u \in\two(D_R)$ that satisfies the equation
\begin{eqnarray}\label{upperbd}
    \begin{array}{rll}
     \fp u + \fq u &\leq \Tilde{\kp}+\Tilde{\kq} \ &\text{in} \ D_R,
        \\
        u& \leq M \doos \ &\text{in} \ D_{2R},
    \end{array}
\end{eqnarray}
for $\Tilde{\kt}>0$ and $M>\widetilde{M}_0>0$ where $t=p,q.$
\begin{proposition}\label{uppbound}
    Let $u \in\two(D_R)$ solve \eqref{upperbd} for $M>\widetilde{M}_0>0$. There exist $\sigma_4\in (0,1],$ $C_6' >1$ depending on $N,p,q,s,\Omega$ and for all $\varepsilon>0, $ a constant $\Tilde{C_\varepsilon}'=\Tilde{C_\varepsilon}'(N,p,q,s,\Omega,\varepsilon)$ such that for all $0<R< \rho/4,$
 we have 
 \begin{equation*}
     \begin{aligned}
         \inf_{D_{R/4}} \left(M-\frac{u}{\doos}\right) \geq& \sigma_4 \Ex- \var \left\|M-\frac{u}{\doos}\right\|_{L^{\infty}(D_R)}-C_6' \tail_1\left(\left(\frac{u}{\doos}-M\right)_{+},2R \right)R^s
         \\
         &- \sum_{t=p,q} \Tilde{C_{\var}}' \left[M+\Tilde{K_t}^{\frac{1}{t-1}}+\tail_{t-1}\left(\left(\frac{u}{\doos}-M\right)_{+},2R\right)\right]R^{\frac{s}{t-1}}.
     \end{aligned}
 \end{equation*}
\end{proposition}
\begin{sketch}
Without loss of generality we assume $M-\frac{u}{\doos}\in L^{\infty}(D_R).$ Fix $\varepsilon>0$ and set $v= u\wedge M\doos. $  Using Proposition \ref{maxminlem} of Appendix and if we write  $\var^{t-1}$ in the place of $\var$ then we get 
       \begin{eqnarray*}
           \begin{array}{rll}
                      \fp v + \fq v &\leq \sum_{t=p,q} K_t + M^{t-2} H_t &\text{ in } D_R,
            \\
            v &\leq M\doos &\text{ in } \RR^N,
           \end{array}
       \end{eqnarray*}
    where 
      \begin{equation*}
        \begin{aligned}
            K_t&=\Tilde{K_t}+\frac{\var^{t-1}}{R^s}\left\|M-\frac{u}{\doos}\right\|_{L^{\infty}(D_R)}^{t-1}+ C'_{\var,t} \tail_{t-1} \left(\left(\frac{u}{\doos}-M\right)_{+},2R \right)^{t-1},
            \\
            H_t&=C'_{2,t} \tail_1\left(\left(\frac{u}{\doos}-M\right)_{+},2R \right).
        \end{aligned}
    \end{equation*}
 We choose 
    \begin{equation*}
        \begin{aligned}
           \sigma_4&= \min\{\sigma_3, \sigma_{\Theta_1}\}<1 ,
           \\
           C&=\max_{t=p,q} \{C_{5,t},C_{\Theta_1,t}\} \geq 1,
        \end{aligned}
    \end{equation*}
    where $0<\sigma_{\Theta_1}\leq 1 \leq C_{\Theta_1,t}$ and $0<\sigma_3\leq 1\leq C_{5,t},$ are given in Lemma \ref{ubgex} and Lemma \ref{usmlex} respectively. We deduce the result through a computation analogous to that in Proposition \ref{lowebound}, thus skipping the detailed proof.
    
\end{sketch}

\section{Weighted H\"{o}lder Regularity For Bounded Non-negative Data}\label{whrp}

In this section, we prove Theorem \ref{mainthmst} when $f$ is a non-negative bounded function and $0<s<1.$ We commence our proof by deriving an estimation of the oscillation of $\frac{u}{\doos}$ near the boundary, when  $u\in \wo(\Om)$ satisfies
\begin{eqnarray}\label{maineqnp}
    \left.
    \begin{array}{rll}
        \fp u + \fq u &= f(x) &\text{ in } \Om,
        \\
       u &=0 &\text{ in } \Om^c,
    \end{array}
    \right\}
\end{eqnarray}
For the given $f,$ we have $0\leq f(x)\leq K $ for some constant $K>0.$ We know that $\frac{u}{\doos}\in L^\infty(\Omega)$ from \cite[Remark 3.7]{GKS23}. Moreover, by \cite[Proposition 2.12]{GKS23},  $\inf_{x\in \Omega}\frac{u}{\doos}>0.$ Using Proposition \ref{lowebound} and Proposition \ref{uppbound} the estimation of the oscillation can be accomplished by modifying the computations outlined in \cite[Theorem 5.1]{IMS20}. Nevertheless we shall provide the outline of the proof in the next Lemma.
\begin{lemma}(Oscillation Lemma) \label{osclemma}
Then there exist $\alpha_1 \in( 0, s],$ $ R_0 \in( 0, \rho / 4)$  and the constant $C_7> \max\{K^\frac{1}{p-1}, K^\frac{1}{q-1}\}$ all depending on $N, p,q, s$ and $\Omega$ such that for all $r \in( 0, R_0)$ ,  
   \begin{equation*}
       \begin{aligned}
          \underset{{D_r(x_1)}}{\osc}  \frac{u}{\doos} \leq C_7  r^{\alpha_1}.
       \end{aligned}
   \end{equation*}
\end{lemma}
\begin{sketch}
     Without loss of generality assume that $x_1 = 0$ and set $v=\frac{u}{\doos}\in L^\infty(\Om),$ $R_0=\min\{1,\rho/4\}.$ Let us define $\widetilde{M}_0=m_0:= \underset{x \in \Omega}{\inf}\frac{u(x)}{\doos(x)}.$ Since $f\geq 0$  by Hopf's lemma \cite[Proposition 2.12]{GKS23}, $m_0>0.$ Now for $R_n=\frac{R_0}{8^n},\;D_n=D_{R_n},\ \widetilde{B}_n= \widetilde{B}_{\frac{R_n}{2}},$ we claim that, there exists $\alpha_1\in (0,s]$ and $\mu\geq 1$ and a nondecreasing sequence $\{m_n\}$ and a nonincreasing sequence $\{M_n\}$ in $\RR$ [all depending on $N,p,q,s,\Omega$] such that 
    \begin{equation}\label{asm}
        \begin{aligned}
            &0<m_0\leq m_n \leqslant \inf_{D_n} v \leqslant \sup_{D_n} v \leqslant M_n \hspace{0.5cm} \text{and}\hspace{0.5cm} M_n - m_n = \mu R_n^{\alpha_1}.
        \end{aligned}
    \end{equation}
We prove our claim by strong induction. Indeed for $n=0,$     
observe by our choice of $m_0,$
     \begin{equation}
         \begin{aligned}
             \inf_{D_0} v > m_0.
         \end{aligned}
     \end{equation}
Set $\|v\|_{L^{\infty}(\Omega)}\leq C_\Omega,$ $M_0:= C_\Omega,$ $\alpha_1\in (0,s]$(to be determined later) and $\mu=\frac{C_\Omega-m_0}{R_{0}^{\alpha_1}}\geq 1$ to verify the first step of induction hypothesis.  Next we assume that our induction hypothesis is satisfied for $n$th step. i.e. 
\begin{equation}\nonumber
        \begin{aligned}
            &0<m_0\ldots \leq m_n \leq M_n\ldots \leq M_0 \hspace{0.5cm} \text{and}\hspace{0.5cm} M_n - m_n = \mu R_n^{\alpha_1}.
        \end{aligned}
    \end{equation}
 If we set $R={R_n}/2$ then $D_{n+1}= D_{R/4},$ and $\widetilde{B}_n=\widetilde{B}_R$, we aim to applying our main estimates in Sections \ref{upbdsec} and \ref{lowbdsec} for $v.$ We note that in the special case when $f\geq 0$ both $m_n , M_n \in [m_0,M_0]$ and $m_0>0.$ Now to complete the proof of oscillation lemma, we need to consider only the case (a) of Theorem 5.1 of \cite{IMS20}. Furthermore, given that $m_n>m_0>0$ and  $M_n>\widetilde{M}_0,$ Propositions \ref{lowebound} and \ref{uppbound} verifies $(5.3)$ and $(5.4)$ in Theorem 5.1 of \cite{IMS20}. Now, proceeding analogously we prove the oscillation lemma. 
\end{sketch}
To establish the final conclusion of our main theorem, we need to invoke a result given in \cite{IMS20}, (see \cite[Theorem 1.2]{RS14} for a proof). Here, we provide the statement of the lemma:
\begin{lemma} \cite[Lemma 5.2]{IMS20}\label{roslem}
     Let $\partial \Omega$ be $C^{1,1}$. If $v \in L^{\infty}(\Omega)$ satisfies the following conditions: 
     \begin{itemize}
         \item[(i)]  $\|v\|_{L^{\infty}(\Omega)} \leqslant C ;$
         \item[(ii)] for all $x_1 \in \partial \Omega, r>0$ we have $\underset{D_r\left(x_1\right)}{\mathrm{osc}} v \leqslant C r^{\beta_1}$;
         \item[(iii)]\, if $\mathrm{d}_{\Omega}\left(x_0\right)=R$, then $v \in C^{\beta_2}\left(B_{R / 2}\left(x_0\right)\right)$ with 
         \begin{equation*}
             [v]_{C^{\beta_2}\left(B_{R / 2}\left(x_0\right)\right)} \leqslant C\left(1+R^{-\mu}\right),
         \end{equation*}
     \end{itemize}
    for some $C_8, \mu>0$ and $\beta_1, \beta_2 \in( 0,1)$, then there exist $\alpha \in( 0,1), C_9>0$ depending on the parameters and $\Omega$ such that $v \in C^\alpha(\bar{\Omega})$ and $[v]_{C^{\alpha}(\bar{\Omega})} \leqslant C_9$.
\end{lemma}
\underline{\textbf{Proof of Theorem \ref{mainthmst} for $(a)$}}  \par
Set $v= \frac{u}{\doos}.$ By \cite[Remark 12]{GKS22}, we get that 
\begin{equation}\label{vprp1}
    \begin{aligned}
        \|v\|_{L^\infty} <C,
    \end{aligned}
\end{equation}
which verifies condition $(i)$ of previous Lemma. 
Using Lemma \ref{osclemma}, we find $C>0$ such that 
\begin{equation}\label{vprp2}
    \begin{aligned}
        \underset{D_r\left(x_1\right)}{\mathrm{osc}} v \leqslant C r^{\alpha_1} \text{ for all } r>0,
    \end{aligned}
\end{equation}
thus verifying condition $(ii)$ of previous Lemma.
Using \cite[Remark 12]{GKS22}, we can guarantee the existence $\alpha_2 \in (0,s]$ such that $u \in C^{\alpha_2}(\Bar{\Om})$ and $\|u\|_{C^{\alpha_2}(\Bar{\Omega})}\leq C,$ for some $C>0.$ Thanks to \cite[Theorem 2.10]{GKS20} and \cite[Theorem 2.6]{GKS23}, we have the following interior regularity result :
\begin{equation}\label{intreg}
    \begin{aligned}
       \relax [u]_{C^{\alpha_2}({B_{R/2}}(x_0))} \leq \frac{C}{R^{\alpha_2}}\left[ (KR^{qs})^\frac{1}{q-1}+ \|u\|_{L^{\infty}(\RR^N)}+ T_{q-1}(u,0,R) + R^\frac{qs-ps}{q-1}\right] ,
    \end{aligned}
\end{equation}
where $B_{R}\subset\subset \Om,$ for all  $x_0 \in \Om,$ $R=\doo(x_0)$ and \begin{equation*}
    \begin{aligned}
       T_{t-1}(u,0,R):&=\left(R^{ts} \int_{B^{c}_{R}} \frac{|u(y)|^{t-1}}{|y|^{N+st}}\right)^\frac{1}{t-1} \text{ for any } t>1.
    \end{aligned}
\end{equation*} 
For any $t>1$ we get that
\begin{equation*}
    \begin{aligned}
       T_{t-1}(u,0,R)\leq
         \|u\|_{L^{\infty}(\RR^N)} R^{t^{\prime}s}\left( \int_{B^{c}_{R}} \frac{1}{|y|^{N+st}}\right)^\frac{1}{t-1} \leq C  \|u\|_{L^{\infty}(\RR^N)}.
    \end{aligned}
\end{equation*} 
Plugging this into \eqref{intreg}, we get 
\begin{equation}\label{uds1}
    \begin{aligned}
      \relax [u]_{C^{\alpha_2}({B_{R/2}}(x_0))}\leq \frac{C}{R^{\alpha_2}}\left[ (KR^{qs})^\frac{1}{q-1}+ \|u\|_{L^{\infty}(\RR^N)}+ C  \|u\|_{L^{\infty}(\RR^N)} + R^\frac{qs-ps}{q-1}\right]  &\leq \frac{C}{R^{\alpha_2}}  ,
    \end{aligned}
\end{equation}
for all $x_0 \in \Om$ and $R=\doo(x_0) \leqslant \operatorname{diam}(\Om).$
For the same choice of $x_0$ and $R$ we get 
\begin{equation}\label{uds2}
    \begin{aligned}
          \relax [\dooso]_{C^{\alpha_2}({B_{R/2}}(x_0))} &\leq \frac{C}{R^{s+\alpha_2}},
    \end{aligned}
\end{equation}
from \cite[pp. 292]{RS14}. Now given that \eqref{uds1} and \eqref{uds2} are proved, we get the following for all $x,y \in B_{R/2}(x_0)$
\begin{equation}\label{vprp3}
    \begin{aligned}
           \frac{|v(x)-v(y)|}{|x-y|^{\alpha_2}} &\leqslant \frac{|u(x)\dooso(x)-u(y)\dooso(x)|}{|x-y|^{\alpha_2}}+ \frac{|u(y)\dooso(x)-u(y)\dooso(y)|}{|x-y|^{\alpha_2}}
           \\ & \leqslant \relax [u]_{C^{\alpha_2}({B_{R/2}}(x_0))} \|\dooso\|_{L^{\infty}(B_{R/2}(x_0))} + \|u\|_{L^{\infty}(\Om)}[\dooso]_{C^{\alpha_2}({B_{R/2}}(x_0))}
           \\ & \leqslant \frac{C}{R^{\alpha_2}} \left(\frac{2}{R}\right)^{s} + \frac{C}{R^{s+\alpha_2}}
          \\ &\leqslant\frac{C}{R^{s+\alpha_2}},
    \end{aligned}
\end{equation}
for chosen $x_0 \in \Om$ and $R=\doo(x_0).$ Thanks to \eqref{vprp1}, \eqref{vprp2} and \eqref{vprp3} all the assumptions of Lemma \ref{roslem} are satisfied with  $\beta_1= \alpha_1,$ $\beta_2= \alpha_2$ and $\mu= s+\alpha_2$ and hence $\|v\|_{C^\alpha(\Bar{\Om})}\leqslant C$ for $\alpha\in (0,s]$ and $C$ depending on  $N, p,q, s$, $\Omega$, $\|f\|_{L^{\infty}(\Omega)}.$ \hfill\qed

\section{Weighted H\"{o}lder Regularity For Bounded sign changing Data}\label{whrs}
In this section, we aim to address the fine boundary regularity for the sign changing bounded data when $s\in (0,\frac{1}{q}).$ Observe that, when the data $f$ is non-negative, oscillation Lemma \ref{osclemma} is proved using Proposition \ref{lowebound} of Section \ref{lowbdsec} and  Proposition \ref{uppbound} of Section \ref{upbdsec}. This relies on the crucial assumption $m_n > m_0 > 0,$ $M_n>\widetilde{M}_0>0$ which is reasonable for non-negative data due to Hopf's Lemma. But when the data $f$ is sign changing, the constant $M$ and $m$ involved in the calculations in Section \ref{lowbdsec} and \ref{upbdsec} cannot have a positive lower bound. To overcome this issue for the sign changing $f$, it is imperative  to understand the exact parameter estimates of solutions to 
\begin{eqnarray}\label{eqnboth}
\begin{array}{rll}
 \fp \vm + \fq \vm&=\mu &\text{ in } U,
            \\
            \vm&=0 &\text{ in } U^{c},
	\end{array} 
\end{eqnarray}
in the case when $\mu\rightarrow 0$ as well as $\mu\rightarrow \infty.$
As discussed in Section 3, using the  transformation $v_\mu\rightarrow \mu^{\frac{-1}{q-1}} v_\mu,$ the Dirichlet problem \eqref{eqnboth} may be converted to 
\begin{eqnarray*}
    \begin{array}{rll}
     \beta  \fp u_\beta + \fq u_\beta &= 1 & \text{ in } U,\\
       u_\beta&=0 & \text{ in } U^c,   \end{array}
\end{eqnarray*}
where $\beta=\mu^{\frac{p-q}{p-1}}.$
Now, the exact behaviour of solutions $v_\mu$ as $\mu\rightarrow \infty$ may be established by using the $C^s_{\loc}$ regularity result given in Appendix \ref{csreg}. This particular transformation does not yield any meaningful result when $\mu\rightarrow 0.$ Consequently, when $\mu\rightarrow 0,$ we consider the transformation $\vm\rightarrow \mu^{\frac{-1}{p-1}} v_\mu,$ which leads to  
\begin{eqnarray}\label{transformedeqn}
    \begin{array}{rll}
        \fp \vb +\beta  \fq \vb &= 1 &\text{ in } U,
      \\  \vb&=0 &\text{ in } U^c,
    \end{array}
\end{eqnarray}
where $\beta\rightarrow 0$ when $\mu\rightarrow 0.$ Unfortunately an interior $C^s$ regularity has not been proved yet for \eqref{transformedeqn}, to the best of our knowledge. Nevertheless, by the Theorem \ref{careg} of Appendix we can obtain a uniform $L^{\infty}$ bound for \eqref{transformedeqn}, thus necessitating consideration of the parameter range $s \in (0,\frac{1}{q})$. Under this additional assumption on $s$ we shall now modify the calculations in Section \ref{lowbdsec} and Section \ref{upbdsec} to establish weighted H\"{o}lder regularity when $s\in (0,\frac{1}{q})$ and $f$ is sign changing. Our first step is to obtain a result analogous to 
Proposition \ref{lowebound} of Section \ref{lowbdsec}. 
In Section \ref{lowbdsec}, we established Proposition \ref{lowebound} under the assumption \( m > m_0 > 0 \) for a fixed \( m_0 > 0 \). In the present section, we do not impose this assumption. Instead, leveraging Lemma \ref{muuniformlemma}, an argument analogous to that in Section \ref{lowbdsec} (where Lemma \ref{muuniformlemma} is employed in place of Lemma \ref{lowbdlem}) yields Proposition \ref{lowebound} under the condition \( s \in (0, \frac{1}{q}) \).
We summarize the result in the following proposition and omit the proof.
\begin{proposition}\label{lowebounde}
Let $s\in (0,\frac{1}{q})$ and $u \in \two(D_R)$ solves 
\begin{eqnarray}\label{localeqn1e}
    \begin{array}{rll}
         \fp u + \fq u &\geq -\Tilde{\kp}-\Tilde{\kq} \ &\text{in} \ D_R,
        \\
        u& \geq m \doos \ &\text{in} \ D_{2R}.
    \end{array}
\end{eqnarray}
for $\Tilde{\kp}, \Tilde{\kq},m\geq 0.$ There exist $0<\sigma_2\leq 1,$  $C_6 >1$ depending on $N,p,q,s,\Omega$ and for all $\varepsilon>0, $ a constant $\Tilde{C_\varepsilon}=\Tilde{C_\varepsilon}(N,p,q,s,\Omega,\varepsilon)>0$ such that for all $0<R< \rho/4,$
 we have 
 \begin{equation*}
     \begin{aligned}
         \inf_{D_{R/2}} (\frac{u}{\doos}-m) \geq& \sigma_2 \ex- \varepsilon \Big\|\frac{u}{\doos}-m \Big\|_{L^{\infty}(D_R)}-C_6 \tail_1\Big(\big(-\frac{u}{\doos}+m\big)_{+},2R \Big)R^s
         \\
         &- \sum_{t=p,q} \Tilde{C_\varepsilon} \Bigg[m+\Tilde{K_t}^{\frac{1}{t-1}}+\tail_{t-1}\Big(\big(-\frac{u}{\doos}+m\big)_{+},2R \Big)\Bigg]R^{\frac{s}{t-1}}.
     \end{aligned}
 \end{equation*}
\end{proposition}
We now modify the arguments of section \ref{upbdsec} to get a result analogous to Proposition \ref{uppbound}. We start with modifying Lemma \ref{lemdoos}.
\begin{lemma} \label{lemdoosm}
   Let $s\in (0,\frac{1}{q})$ and $E_R$ is defined in Section \ref{upbdsec}. Given an $M\in(0,m_0)$ there exists a $\Lambda>1$ large enough such that if any function $\uM\in \wo(E_R)$ satisfy
    \begin{eqnarray*}
        \begin{array}{rll}
       \fp\uM+\fq\uM&=\frac{\Lambda(\Mpq)}{R^s} &\text{ in } E_R,
            \\
            \uM&=0 &\text{ in } E^{c}_R,
        \end{array}
    \end{eqnarray*}
  
    implies that 
    \begin{equation*}
        \begin{aligned}
            u_M(x)\geq M \doos \text{ for all } x\in D_{3R}\setminus D_R.
        \end{aligned}
    \end{equation*}
\end{lemma}
 \begin{proof}
Let $\widetilde{v}_M(x)=\frac{\uM(Rx)}{M R^s}$ for $x\in E_1.$ Then
\begin{eqnarray*}\label{veqnmum}
  \left\{  \begin{array}{rll}
        \fp \widetilde{v}_M + M^{q-p}\fq \widetilde{v}_M &= \Lambda(1+M^{q-p}) &\text{ in }E_1,
         \\
         \widetilde{v}_M&=0 &\text{ in }E_{1}^{c}.
    \end{array}\right.
\end{eqnarray*}
    Also, let $\vM\in \wo(E_1)$ solves 
\begin{eqnarray*}
  \left\{  \begin{array}{rll}
          \fp \vM + M^{q-p}\fq \vM &=  \Lambda &\text{ in }E_{1},
         \\
         \vM&=0 &\text{ in }E_{1}^{c}.
    \end{array}\right.
\end{eqnarray*} 
Using $\vM$ as a test function in the weak formulation , and the continuous embedding $\wsp \hookrightarrow L^1$ for any $p>1,$ we get $[\vM]_{W_{0}^{s,p}}\leq C,$ independent of $M.$
Hence $\vM\rightharpoonup v_0$ in $\wsp(E_1)$ upto a subsequence. Using Theorem \ref{careg} in Appendix, we have that $\|\vM\|_{L^{\infty}_{\loc}(E_1)}$ is uniformly bounded for small values of $M.$ Hence, for $s\in (0,\frac{1}{q}),$ we can show that 
    \begin{equation*}
        \begin{aligned}
         |   \int_{\mathbb{R}^N \times \mathbb{R}^N } \dfrac{(\vM(x)-\vM(y))^{q-1}(\varphi(x)-\varphi(y))}{|x-y|^{N+sq}}  \,dx\,dy |\leq C 
        \end{aligned}
    \end{equation*}
    for each $\varphi \in C_c^\infty(\Om)$ where $C$ is independent of $M.$ Clearly as $M\rightarrow 0,$ 
     \begin{equation*}
        \begin{aligned}
            M^{q-p}\int_{\mathbb{R}^N \times \mathbb{R}^N } \dfrac{(\vM(x)-\vM(y))^{q-1}(\varphi(x)-\varphi(y))}{|x-y|^{N+sq}}  \,dx\,dy \rightarrow 0 \mbox{ for all } \varphi\in C_c^\infty(E_1).
        \end{aligned}
    \end{equation*}
    Due to the weak-weak continuity property \cite[Lemma 2.2]{CMM18} of fractional $p$-Laplacian and the density argument,  the function $v_0$ solves
    \begin{equation}
        \begin{aligned}
            \fp v_0 &=\Lambda \ \text{in} \ E_{1}\\
v_0 &= 0 \hspace{0.2cm} \text{in} \ E_{1}^{c}.
        \end{aligned}
    \end{equation}
    Thanks to the $L^{\infty}_{\loc}(\Om)$ uniform bound, we can use Ascoli-Arzela theorem to infer that $\vM$ has a uniformly convergent subsequence. Uniqueness of the solution of \eqref{evo} would imply that   $\|\vM - v_0\|_{L^\infty(K_1)}\rightarrow0$  as $M \to 0$ for any fix compact subset $K_1$ of $E_1.$ Now, 
   using the strong comparison principle for fractional $p$-Laplacian,  we have  $\underset{K_1}{\inf} \, v_0\geq C_1 \Lambda^{\frac{1}{p-1}} >0.$ Since $\vM$ uniformly converges to $v_0,$ for sufficiently small $M$, we can conclude that 
   \begin{equation*}
    \begin{aligned}
      \vM(x)\geq v_0(x)-\varepsilon \geq C_1 \Lambda^{\frac{1}{p-1}}-\varepsilon \ \text{for all } \ x\in  K_1.
    \end{aligned}
\end{equation*}
 Now, choosing $\varepsilon$ small enough we obtain, 
\begin{equation}\label{vmeqn}
    \begin{aligned}
        \underset{K_1}{\inf} \, \vM \geq C\Lambda^{\frac{1}{p-1}}>0  \, \text{  for all } M< m_0 \,\, \text{ and } \forall x\in K_1 ,
    \end{aligned}
\end{equation}
where $C=C(N,s,p,\Omega)$ is a constant independent of $M.$ Now, let  
\begin{equation*}
    \begin{aligned}
        k=\min_{t=p,q}\left\{\left(\frac{\Lambda}{2\max_{E_1\setminus K_1}( \ft \deones)}\right)^\frac{1}{t-1}, \frac{C \Lambda^{\frac{1}{p-1}}}{\max_{K_1}\deones} \right\}.
    \end{aligned}
\end{equation*}
 Now using comparison principle we obtain
 \begin{equation}\label{deonevm}
     \begin{aligned}
         k \deones \leq \vM \text{ in } \RR^N.
     \end{aligned}
 \end{equation}  
 We note that $v_M\leq \widetilde{v}_M$ in $\RR^N$ where $\widetilde{v}_M(x)=\frac{\uM(Rx)}{M R^s}$ for $x\in E_1.$ Thus, 
\begin{equation} \label{umdoos}
    \begin{aligned}
        \frac{\uM(x)}{M }\geq k \ders \geq k \, C_5 \doos \text{ in } D_{3R}\setminus D_R.
    \end{aligned}
\end{equation}
We can choose $\Lambda$ large enough such that $k \, C_5  > 1 $ and hence the result.

 \end{proof}
 Once we proved Lemma \ref{lemdoosm}, we can modify Lemma \ref{barlem1} to construct the  barrier for $M\in(0,m_0).$
 \begin{lemma}
  Let $s\in (0,\frac{1}{q}),$ $\Bar{x}\in D_{R/2},$   and $R\in(0,\rho/4).$ Given an $M\in (0,m_0),$ there exists a function $v\in \wo(\Om)\cap C(\RR^N)$ and a positive constant $C(N,p,q,s,\Omega)>1$ such that the following conditions are satisfied: 
    \begin{itemize}
        \item[(i)] $\fp v + \fq v\geq -\frac{C(M^{p-1}+M^{q-1})}{R^s}$ in $D_{2R}$;
         \item[(ii)] $v(\Bar{x})=0;$ 
          \item[(iii)] $v\geq M\doos$ in $D^{c}_R;$
           \item[(iv)] $|v|\leq C R^s $ in $D_{2R}.$
    \end{itemize}
\end{lemma}
\begin{sketch}
First, we construct the lower obstacle.  Define $\va \in \wo(E_R)$ such that 
    \begin{eqnarray*}
        \begin{array}{rll}
           \fp \va+ \fq \va&= \frac{\Lambda(\Mpq)}{R^s} \ &\text{in} \ E_R,\\
\va &= 0 \ &\text{in} \ E_{R}^{c} .
        \end{array}
    \end{eqnarray*}  
     where $\Lambda>0.$ Using Lemma \ref{lemdoosm}, and an analogue to Lemma \ref{lowbdlem}, for the given $M>0$  we  choose $\Lambda$ large enough to get 
      \begin{equation} \label{lowobsp31}
        \begin{aligned}
            \va \geq  M\doos \text{ in }  D_{3R}\setminus D_R, \text{ and } \va\leq C_\Lambda M R^s \text{ in } \RR^N.
        \end{aligned}
    \end{equation}

    This $\varphi$ would serve us as the lower obstacle. Fixing $\Lambda>0$ as in \eqref{lowobsp31} we will now construct the upper obstacle. 
  Let $\Psi \in \wo(B_{R/8})$ satisfies
    \begin{eqnarray*}
        \begin{array}{rll}
        \fp \Psi + \fq \Psi&= \frac{\widetilde{\Lambda}(\Mpq)}{R^s} \ &\text{in} \  B_{R/8},\\
\Psi &= 0 \hspace{2cm} &\text{in} \ B_{R/8}^{c}. 
        \end{array}
    \end{eqnarray*}
 We define
\begin{equation*}
\begin{aligned}
        \psi(x)= \max_{\RR^N} \Psi- \Psi(x-\x) \text{ for }x\in\RR^N. 
\end{aligned}
\end{equation*}
  Clearly $\psi \in \tw(\Om),$ $\psi\geq 0,$ and $\psi(\x)=0.$ Moreover we want to prove 
  \begin{equation*}
    \begin{aligned}
            \va \leq C_\Lambda M R^s\leq \psi \text{ in } \RR^N \text{ for } \widetilde{\Lambda} \text{ large enough.}
    \end{aligned}
\end{equation*}
Observing $ \va(x)=0$ in $D_{3R/4},$ we only consider $D^{c}_{3R/4},$ where $\psi(x)=\max_{\RR^N} \Psi.$ 
Similar to \eqref{deonevm} and the first inequality of \eqref{umdoos} of Lemma \ref{lemdoosm}, we have $\widetilde{k} \dis_{B_{\frac{R}{8}}}^{s} (x)\leq \frac{ \Psi(x)}{M}$ for  $x\in \RR^N,$ where 
\begin{equation*}
    \begin{aligned}
        \widetilde{k}:= \min_{t=p,q}\left\{\left(\frac{\widetilde{\Lambda}}{2\max_{B_{\frac{1}{8}}\setminus K_1}( \ft \dis_{B_{\frac{1}{8}}}^{s})}\right)^\frac{1}{t-1}, \frac{C \widetilde{\Lambda}^{\frac{1}{p-1}}}{\max_{K_1}\dis_{B_{\frac{1}{8}}}^{s}} \right\} \text{ for a given } K_1\subset\subset B_{\frac{1}{8}}.
    \end{aligned}
\end{equation*}
So we can choose $\widetilde{\Lambda}$ large enough such that $8^{s} C_\Lambda \leq \widetilde{k}$ to conclude that
\begin{equation*}
    \Psi(x) \geq 8^{s} C_\Lambda M \dis_{B_{\frac{R}{8}}}^{s} (x)\text{ for } x\in \RR^N,
\end{equation*}
for $\widetilde{\Lambda}$ large enough and henceforth fix $\widetilde{\Lambda}.$  Observe that we now constructed our both obstacles $\phi$ and $\psi$ and consequently define $\Phi\in \wo(\Omega)$ as a unique minimizer  
\begin{equation*}
    \begin{aligned}
      \Phi:=\min \left\{\frac{1}{p}[u]^{p}_{s,p}+\frac{1}{q}[u]^{q}_{s,q}: u\in \wo(\Omega)\text{ and } \va \leq u \leq \psi \text{ in } \RR^N\right\}  ,
    \end{aligned}
\end{equation*}
which satisfies the properties (a)-(d) as described in Lemma \ref{barlem1}. Now we can modify the calculations of Lemma \ref{barlem1} to construct a $v$ satisfying our hypothesis.
\end{sketch}
Hence we can prove that an analogue to Proposition \ref{uppbound} for $s\in  (0,\frac{1}{q})$ and $M\in (0,m_0).$
\begin{proposition}\label{uppbounds}
    Let $s\in (0,\frac{1}{q})$ and  $u \in\two(D_R)$ solve 
    \begin{eqnarray}
    \begin{array}{rll}
     \fp u + \fq u &\leq \Tilde{\kp}+\Tilde{\kq} \ &\text{in} \ D_R,
        \\
        u& \leq M \doos \ &\text{in} \ D_{2R},
    \end{array}
\end{eqnarray}
for $\Tilde{\kp},\Tilde{\kq},M > 0.$ There exist $\sigma_4\in (0,1],$ $C_6' >1$ depending on $N,p,q,s,\Omega$ and for all $\varepsilon>0, $ a constant $\Tilde{C_\varepsilon}'=\Tilde{C_\varepsilon}'(N,p,q,s,\Omega,\varepsilon)$ such that for all $0<R< \rho/4,$
 we have 
 \begin{equation*}
     \begin{aligned}
         \inf_{D_{R/4}} \left(M-\frac{u}{\doos}\right) \geq& \sigma_4 \Ex- \var \left\|M-\frac{u}{\doos}\right\|_{L^{\infty}(D_R)}-C_6' \tail_1\left(\left(\frac{u}{\doos}-M\right)_{+},2R \right)R^s
         \\
         &- \sum_{t=p,q} \Tilde{C_{\var}}' \left[M+\Tilde{K_t}^{\frac{1}{t-1}}+\tail_{t-1}\left(\left(\frac{u}{\doos}-M\right)_{+},2R\right)\right]R^{\frac{s}{t-1}}.
     \end{aligned}
 \end{equation*}
\end{proposition}
Now, our attention turns to the completing proof of the second part Theorem \ref{mainthmst}, which demonstrates the weighted H\"{o}lder continuity properties of the solutions of problem \eqref{eqn1}. In accordance with a customary strategy, we commence by deriving an estimation of the oscillation of $\frac{u}{\doos}$ near the boundary, where $u\in \wo(\Om)$ satisfies,
\begin{eqnarray}\label{maineqn}
    \left.
    \begin{array}{rll}
        -K\leqslant \fp u + \fq u &\leqslant K &\text{ in } \Om,
        \\
       u &=0 &\text{ in } \Om^c,
    \end{array}
    \right\}
\end{eqnarray}
with some $K>0.$ Given that Proposition \ref{lowebounde} and Proposition \ref{uppbounds} have been previously established, it becomes evident that the estimation of the oscillation can be accomplished by modifying the computations outlined in \cite[Theorem 5.1]{IMS20}. Consequently, we omit the detailed proof here.
\begin{lemma}(Oscillation Lemma) \label{osclemmas}
       Let $s\in (0,\frac{1}{q}),$  $ x_1 \in \partial \Omega$ and  $u \in\wo(\Omega)$ solve \eqref{maineqn}. There exist $\alpha_1 \in( 0, s],$ $ R_0 \in( 0, \rho / 4)$  and $C_7> \max\{K^\frac{1}{p-1}, K^\frac{1}{q-1}\}>1$ all depending on $N, p,q, s$ and $\Omega$ such that for all $r \in( 0, R_0)$ ,  
   \begin{equation*}
       \begin{aligned}
          \underset{{D_r(x_1)}}{\osc}  \frac{u}{\doos} \leq C_7  r^{\alo}.
       \end{aligned}
   \end{equation*}
\end{lemma}
\underline{\textbf{Proof of Theorem \ref{mainthmst} for $(b)$ }}
\\
The proof methodology is based on that of \cite[Theorem 1.1]{IMS20}, with the detailed steps modified analogously to Theorem \ref{mainthmst} under condition $(a)$.  First we set $v= \frac{u}{\doos}.$ By \cite[Remark 12]{GKS22}, we get that $  \|v\|_{\infty} <C.$ Once Lemma \ref{osclemmas} is proved, we can find $C>0$ such that $\underset{D_r\left(x_1\right)}{\mathrm{osc}} v \leqslant C r^{\alpha_1} \text{ for all } r>0.$ Similar to \eqref{uds1}
we get $\relax [u]_{C^{\alpha_2}({B_{R/2}}(x_0))}\leq \frac{C}{R^{\alpha_2}}.$
Now using Lemma \ref{roslem} we get our result. \hfill\qed
\begin{remark}
Let $\beta\in(0,\beta_{0})$ and $\vb$ solves
\begin{eqnarray*}
    \begin{array}{rll}
        \fp \vb +\beta  \fq \vb &= 1 &\text{ in } U,
      \\  \vb&=0 &\text{ in } U^c,
    \end{array}
\end{eqnarray*}
for a  bounded smooth domain $U.$ Results presented in Appendix Theorem \ref{careghld} guarantees that $\vb$ is uniformly bounded in $C^{\alpha_0}_{\loc}$ for some $\alpha_0>0.$ Our analysis proves that, if \(\vb \in C^{\alpha}_{\loc}(\Omega)\) uniformly in \(\beta\), then a result analogous to Theorem \ref{mainthmst}, part \((b)\), holds  for the range \(0 < s < \frac{\alpha(q-1)+1}{q}\). In particular, if \(\alpha > \frac{sq - 1}{q - 1}\), the restriction \(0 < s < \frac{1}{q}\) can be removed, allowing any \(s \in (0,1)\).
\end{remark}

\appendix
\counterwithin{theorem}{section}
\counterwithin{proposition}{section}
\counterwithin{lemma}{section}
\counterwithin{equation}{section}

\section{}\label{app}
We shall first state a uniform interior H\"{o}lder regularity result which is found helpful in passing through the limit in Section \ref{est}. 
\begin{theorem}\label{csreg}
Let $\mu_0>0$ be given and $\vmu \in \mathcal{W}_0(\Om)$ be the weak solution to the problem $(P_\mu),$ defined as
\begin{eqnarray*}
    (P_\mu)\left\{
    \begin{array}{rll}
        \mu \fp \vmu + \fq \vmu &= f(x) &\text{ in } \Omega,
      \\  \vmu&=0 &\text{ in } \Omega^c,
    \end{array}
    \right.
\end{eqnarray*}
for $f\in L^{\infty}(\Om)$  and $\mu\in(0, \mu_0).$ Then for every $\sigma \in (0,  s)$ we have $\vmu \in C^{0,\sigma}_{\loc}(\Om)$ and for any given compact subset $K$ of $\Omega,$
\begin{equation*}
   \relax \|\vmu\|_{C^\sigma(K)} \leqslant C(K, N,s,p,q,\sigma,\mu_0)
\end{equation*}
where  $C$ is independent of $\mu.$ In particular, $\vmu\in C^{0,\sigma_{0}}_{\loc}(\Omega)$ for any $\frac{sq-1}{q-1}<\sigma_{0}<s.$
\end{theorem}
\begin{sketch}
    First we shall prove the local boundedness of solution when $f\in L^{\infty}(\Omega)$ and follow the proof of \cite[Lemma 3.2 and Proposition 3.3]{GKS23}. Set $\mu_{l}=\mu$ if $l=p$ and $\mu_{l}=1$ for $l=q.$ We infer the Caccioppoli inequality \cite[ Eqn No. (3.1)]{GKS23}  as 
    \begin{equation}\label{ccp}
        \begin{aligned}
            \int_{E}\int_{E}\frac{|\omega(x)\psi(x)-\omega(y
            )\psi(y)|^{q}}{|x-y|^{N+sq}}\thinspace dx\thinspace dy \leq& C  \sum_{l\in\{p,q\}}  \mu_{l} \int_{E} \int_{E} \frac{|\psi(x)-\psi(y)|^{l}}{|x-y|^{N+sl}}(\omega(x)^{l}+ \omega(y)^{l}) \thinspace dx\thinspace dy
            \\ 
            +& C  \sum_{l\in\{p,q\}} \mu_{l} (\sup_{y\in \supp \psi} \int_{\RR^{N}\setminus E} \frac{\omega(x)^{l-1}\thinspace dx}{|x-y|^{N+sl}}\int_{E} \omega\psi^{q})
            \\
            +& C \int_{\Omega^{\prime\prime}} |f| \omega \psi^{q}
        \end{aligned}
    \end{equation}
    where $E\Subset\Omega^{\prime\prime}\Subset\Omega,$ $\omega=(\vmu-k)_{\underset{-}{+}}$ for $k \in \RR$ and $\psi\in C_{c}^{\infty}(E)$ with $0\leq \psi\leq 1.$  Indeed we take $\omega\psi^{q}$ as a test function and get 
    \begin{equation*}
        \begin{aligned}
          \mu  A_{p}(\vmu,\omega\psi^{q},\RR^{2N}) +  A_{q}(\vmu,\omega\psi^{q},\RR^{2N})=\int_{\Omega} f \omega \psi^{q}.
        \end{aligned}
    \end{equation*}
    Now following the argument similar to \cite[pp.9-10]{GKS23} we get that 
    \begin{equation*}
        \begin{aligned}
            &\sum_{l\in\{p,q\}}  \mu_{l} \int_{E} \int_{E} \frac{(\omega(x)-\omega(y))^{l-1}(\omega(x)\psi(x)^{q}-\omega(y)\psi(y)^{q})}{|x-y|^{N+ls}}dx\thinspace dy
            \\
            &\leq 2 \sum_{l\in\{p,q\}}  \mu_{l} (\sup_{y\in \supp \psi} \int_{\RR^{N}\setminus E} \frac{\omega(x)^{l-1}\thinspace dx}{|x-y|^{N+sl}}\int_{E} \omega\psi^{q})
            + C \int_{\Omega^{\prime\prime}} |f| \omega \psi^{q},
        \end{aligned}
    \end{equation*}
    in the place of \cite[Eqn. no (3.5)]{GKS23}. Now following the arguments of \cite[pp. 10]{GKS23}, instead of \cite[Eqn. no (3.7)]{GKS23}, we get that 
    \begin{equation*}
        \begin{aligned}
            &\int_{E}\int_{E}\frac{|\omega(x)-\omega(y
            )|^{q}}{|x-y|^{N+sq}}(\psi(x)^{q}+\psi(y)^{q})\thinspace dx\thinspace dy \\&\leq C  \sum_{l\in\{p,q\}}  \mu_{l} \int_{E} \int_{E} \frac{|\psi(x)-\psi(y)|^{l}}{|x-y|^{N+sl}}(\omega(x)+ \omega(y))^{l} \thinspace dx\thinspace dy
            \\ 
            &+ C  \sum_{l\in\{p,q\}} \mu_{l} (\sup_{y\in \supp \psi} \int_{\RR^{N}\setminus E} \frac{\omega(x)^{l-1}\thinspace dx}{|x-y|^{N+sl}}\int_{E} \omega\psi^{q})
            \\
            &+ C \int_{\Omega^{\prime\prime}} |f| \omega \psi^{q}.
        \end{aligned}
    \end{equation*}
    Thus we can conclude the \eqref{ccp}. Since $\mu\in(0,\mu_{0}),$ the equation \eqref{ccp} reads as
    \begin{equation*}
        \begin{aligned}
            &\int_{E}\int_{E}\frac{|\omega(x)\psi(x)-\omega(y
            )\psi(y)|^{q}}{|x-y|^{N+sq}}\thinspace dx\thinspace dy \\\leq& C \max\{\mu_{0},1\} \Big( \sum_{l\in\{p,q\}}   \int_{E} \int_{E} \frac{|\psi(x)-\psi(y)|^{l}}{|x-y|^{N+sl}}(\omega(x)^{l}+ \omega(y)^{l}) \thinspace dx\thinspace dy
            \\ 
            +&   \sum_{l\in\{p,q\}}  (\sup_{y\in \supp \psi} \int_{\RR^{N}\setminus E} \frac{\omega(x)^{l-1}\thinspace dx}{|x-y|^{N+sl}}\int_{E} \omega\psi^{q})
            + \int_{\Omega^{\prime\prime}} |f| \omega \psi^{q} \Big).
        \end{aligned}
    \end{equation*}
    Therefore we end up with the \cite[Eqn. no (3.1)]{GKS23}. Similar to \cite[Definition 2.2]{GKS23} set
    \begin{equation*}
        \begin{aligned}
            T_{m,\alpha}(u;x_0,R)&:=\left(R^{\alpha}\int_{B_{R}(x_0)}^{c}\frac{|u(y)|^{m}}{|x_0-y|^{N+\alpha}}dy\right)^{\frac{1}{m}},
            \\
            T_{m-1}(u;x,R)&:=T_{m-1,sm}(u;x,R).
        \end{aligned}
    \end{equation*} For $B_{r}(x_0)\Subset E$ where $x_{0}\in \Omega,$ we get \cite[Proposition 3.3]{GKS23} as 
    \begin{equation}\label{localbd}
        \begin{aligned}
            \sup_{B_{r/2}(x_{0})} \omega \leq  C (\intavg_{B_{r}} \omega^{q})^{\frac{1}{q}} + T_{p-1}(\omega;x_0;\frac{r}{2})^{\frac{p-1}{q-1}} + T_{q-1}(\omega;x_0;\frac{r}{2})+1
        \end{aligned}
    \end{equation}
    where $C= (C(N,p,q,s)(1+\|f\|_{L^{\infty}}))^{\frac{q_{s}^{*}}{q^{2}}}.$ Using $(P_{\mu})$ we have 
    \begin{equation*}
        \begin{aligned}
            \relax [\vmu]_{\wsq}^{q}\leq \int_{\Omega} f \vmu,
        \end{aligned}
    \end{equation*}
    thus, $\|\vmu\|_{\wsq}$ is bounded independent of $\mu.$ Therefore, $T_{q-1}(\vmu;x_0;\frac{r}{2})$ is bounded independent of $\mu.$ Since $p<q,$ and $\vmu=0$ in $\Omega^{c},$ this will imply $T_{p-1}(\vmu;x_0;\frac{r}{2})$ is bounded independent of $\mu.$ For $k=0$ i.e. $\omega=\vmu$ we can bound the right hand side of \eqref{localbd} independent of $\mu$ hence $\|\vmu\|_{L^{\infty}_{\loc}}$ is independent of $\mu.$
    \\
    Next, we shall modify the proof of \cite[Theorem 2.6]{GKS23} to get our desired interior Holder regularity result. In fact, we will show that even under this modified condition, \cite[Proposition 3.9]{GKS23} remains valid. In particular, we shall prove that \cite[eqn. (3.18)]{GKS23} holds independent of $\mu.$ To achieve this, we will reference the proof of \cite[Proposition 3.9]{GKS23} and identify the points of modification. 
    \\
    In our case \cite[eqn. (3.14)]{GKS23} reads as 
    \begin{equation*}
        \begin{aligned}
            &\mu\int_{\RR^{2N}} ([{\vmu}_{h}(x)-\vmu_{h}(y)]^{p-1}-[\vmu(x)-\vmu(y)]^{p-1})(\phi(x)-\phi(y)) \thinspace d\mu_{1}\\
            &+  \int_{\RR^{2N}} ([\vmu_{h}(x)-\vmu_{h}(y)]^{q-1}-[\vmu(x)-\vmu(y)]^{q-1})(\phi(x)-\phi(y)) \thinspace d\mu_{2}
            \\
            &=\int_{B_{2R_{0}}} (f(x+h)-f(x))\phi(x)\thinspace dx
        \end{aligned}
    \end{equation*}
    where $\vmu_h$ is defined in \cite[Proposition 3.9]{GKS23}.
    Following the same argument we shall get the following instead of \cite[eqn. (3.15)]{GKS23}
    \begin{equation*}
        \begin{aligned}
            I_{1}(q,s)\leq - \mu I_{1}(p,s) + |I_{2}(q,s)| + \mu|I_{2}(p,s)|+ |I_{3}(q,s)|+ \mu |I_{3}(p,s)| + |I_{4}(f)|,
        \end{aligned}
    \end{equation*}
    for $I_1,I_2,I_3,I_4$ are defined in \cite[pp. 16]{GKS23}. We can obtain 
    \begin{equation*}
        \begin{aligned}
        \left[\frac{[\delta_{h}\vmu]^{\frac{\beta+q-1}{q}}\eta}{|h|^{\frac{1+\nu \beta}{q}}}\right]_{\wsq(B_{R})}^{q}\leq C\Big[I_{11}(q)+I_{12}(q)+ \mu I_{11}(p)+ \sum_{l\in\{p,q\}}\mu_{l}(|I_{2}(l)|+|I_{3}(l)|+|I_{4}(f)|)\Big]
        \end{aligned}
    \end{equation*}
    in the place of \cite[eqn. (3.18)]{GKS23}. For all previously undefined mathematical notation in the equation, readers are referred to \cite[Proposition 3.9]{GKS23}. Since $\mu\in (0,\mu_{0})$ we get that 
        \begin{equation}\label{pro39scv}
        \begin{aligned}
        \left[\frac{[\delta_{h}\vmu]^{\frac{\beta+q-1}{q}}\eta}{|h|^{\frac{1+\nu \beta}{q}}}\right]_{\wsq(B_{R})}^{q}\leq C_1\Big[I_{11}(q)+I_{12}(q)+  I_{11}(p)+ \sum_{l\in\{p,q\}}(|I_{2}(l)|+|I_{3}(l)|+|I_{4}(f)|)\Big]
        \end{aligned}
    \end{equation}
    where $C_1:= C \max\{\mu_{0},1\}$ is independent of $\mu.$ Observe that \eqref{pro39scv} is exactly the same as the \cite[eqn. (3.18)]{GKS23}. Hence we can follow the arguments provided in \cite[Section 3.1]{GKS23} to conclude that 
    \begin{equation*}
        \begin{aligned}
            \relax[\vmu]_{C^{\sigma}(B_{\frac{R_{0}}{2}}(x_{0}))}\leq (C K_{2}(\vmu))^{i_{\infty}}([\vmu]_{W^{s,q}B_{R_{0}}(x_{0})}+1)^{2-\frac{\sigma}{s}}
        \end{aligned}
    \end{equation*}
    where $i_{\infty}\in\mathbb{N}$ such that $i_{\infty}>\frac{N}{s-\sigma},$ $C(N,s,q,\sigma)>0$ and 
    \begin{equation*}
        \begin{aligned}
            K_{2}(\vmu)= 1+T_{p-1}(\vmu;x_0;\frac{r}{2})^{p-1} + T_{q-1}(\vmu;x_0;\frac{r}{2})^{q-1}+ \|\vmu\|_{L^{\infty}(B_{R_{0}(x_0)})}^{\frac{(q+i_{\infty})(q-1)}{q-2}}+ \|f\|_{L^{\infty}(B_{R_{0}(x_0)})}.
        \end{aligned}
    \end{equation*}
    Since $\|\vmu\|_{\wsq}$ is bounded independent of $\mu$ and $\vmu=0$ in $\Omega^{c}$ we infer that $ K_{2}(\vmu)$ is bounded  independent of $\mu.$ Hence we conclude the result.
\end{sketch}

Next we shall consider a problem very similar to $(P_\mu)$ and state a regularity result which is found useful in Section 7 of this article.
\begin{theorem}\label{careg}
Let $\mu_0>0$ be given and $\zmu \in \mathcal{W}_0(\Om)$ be the weak solution to the problem $(Q_\mu),$ defined as
\begin{eqnarray*}
    (Q_\mu)\left\{
    \begin{array}{rll}
        \fp \zmu + \mu \fq \zmu &= f(x) &\text{ in } \Omega,
      \\  \zmu&=0 &\text{ in } \Omega^c,
    \end{array}
    \right.
\end{eqnarray*}
for $f\in L^{\infty}(\Om), f\geq 0$  and $\mu\in(0, \mu_0).$ The family of functions $\zmu$ is uniformly bounded in $L^\infty(\Omega)$ independent of $\mu.$ 
\end{theorem}
\begin{sketch}
    We follow the idea of \cite[Theorem 2.3]{GKS20} to prove the local boundedness result. Let $r<p_{s}^{*}$ and
    \begin{equation*}
        \begin{aligned}
            \rho_{\mu}&\geq\max\{1,\|\zmu\|_{L^{r}(\Omega)}^{-1}\}, \quad v=(\rho_{\mu}\|\zmu\|_{L^{r}(\Omega)})^{-1}\zmu,
            \\w_{k}&=(v-1+2^{-k})^{+}, \hspace{0.9cm} w_0=v^{+},
        \end{aligned}
    \end{equation*}
    similar to \cite[pp. 7]{GKS20}. Observe that $\|\zmu\|_{\wpp}$ is bounded independent of $\mu,$ and hence $U_{k}=\|w_k\|_{L^{r}}^{r}$ is also bounded independent $\mu.$ Then the \cite[eqn. no (2.2)]{GKS20} becomes 
    \begin{equation*}
        \begin{aligned}
            \|w_{k+1}\|^{p}_{\wpp} \leq C(\rho \|u\|_{L^{r}(\Omega)})^{1-p}\left(\int_{\Omega_{k+1}}f(x)w_{k+1}\right).
        \end{aligned}
    \end{equation*}
    Now following the arguments of \cite[pp. 8-9]{GKS20}( which is actually easier as our right hand side does not depend on $u$) we get $U_{k}\rightarrow 0$ as $k\rightarrow\infty.$ This would imply 
    \begin{equation*}
        \begin{aligned}
            &\|\zmu\|_{L^{\infty }(\Omega)}\leq\rho_{\mu}\|\zmu\|_{L^{r}(\Omega)} \text{ where }
            \\
            & \rho_{\mu}=\max\{1,\|\zmu\|_{L^{r}(\Omega)}^{-1}, (\|\zmu\|_{L^{r}(\Omega)}^{\frac{r^2}{p}-r}\eta^{-1})^{\frac{1}{\gamma}}, C^{\frac{N^{2}}{rs}}\}.
        \end{aligned}
    \end{equation*}
Since $\|\zmu\|_{\wpp}$ is bounded independent of $\mu,$ $\rho_{\mu}\|\zmu\|_{L^{r}(\Omega)}$ is bounded independent of $\mu$ and hence the result.
\end{sketch}
Next prove that for $f(x)\equiv 1$ the solutions of $(Q_{\mu})$ is interior Holder continuous independent of $\mu.$  We follow the approach of \cite[Theorem 1.1]{BOS22}. Indeed we define 
   \begin{equation*}
       \begin{aligned}
           H(x,y,\tau)&=\frac{\tau^{p}}{|x-y|^{sp}}+\mu \frac{\tau^{q}}{|x-y|^{sq}}
           \\
           h(x,y,\tau)&=\frac{\tau^{p-1}}{|x-y|^{sp}}+\mu \frac{\tau^{q-1}}{|x-y|^{sq}}
       \end{aligned}
   \end{equation*}
for $x,y\in\RR^{N}$ and $\tau\geq 0$ similar to \cite[Eqn no. (2.1),(4.1)]{BOS22}.
\begin{lemma}\label{caclem}
    We get that 
   \begin{equation*}
       \begin{aligned}
           &\int_{B_{r}} \int_{B_{r}} H(x,y,|w_{\underset{-}{+}}(x)-w_{\underset{-}{+}}(y)|)(\phi^{q}(x)+\phi^{q}(y))\frac{dx\;dy}{|x-y|^{N}}
           \\
           & \leq  c\int_{B_{r}} \int_{B_{r}} H(x,y,|(w_{\underset{-}{+}}(x)+w_{\underset{-}{+}}(y))(\phi(x)-\phi(y))|)\frac{dx\;dy}{|x-y|^{N}}
           \\
           & + c \left(\sup_{x\in\supp \phi}\int_{\RR^{N}\setminus B_{r}}h(x,y,w_{\underset{-}{+}}(y)\frac{dy}{|x-y|^{N}})\right)\int_{B_{r}}w_{\underset{-}{+}}(x)\phi^{q}(x)\thinspace dx + \int_{B_{2r}} w_{\underset{-}{+}}(x)\phi^{q}(x)\thinspace dx
       \end{aligned}
   \end{equation*}
   for $c(N,s,p,q,\Omega)$ a constant and $B_{2r}\Subset\Omega$ a ball, $\phi\in C_{0}^{\infty}(B_{r})$ with $0\leq \phi \leq 1,$ $w_{\underset{-}{+}}=(\zmu-k)_{\underset{-}{+}}$ where $k\geq 0.$ 
\end{lemma}
\begin{sketch}
    Use $w_{\underset{-}{+}}\phi^{q}$ as the test function and follow the arguments of  \cite[ Lemma 4.2]{BOS22} to conclude this result.
\end{sketch}
Next fix $\Omega^\prime\Subset\Omega$ and define $M\equiv M(\Omega^\prime)= 1+\|u\|_{L^{\infty}(\Omega^\prime)}^{q-p}.$ We start with obtaining a logarithmic type estimate.
\begin{lemma}\label{loglem}
    Let $u\in L_{\loc}^{\infty}(\Omega)$ be a supersolution to $(Q_{\mu})$ and $B_{R}(x_0)\subset\Omega^{\prime}\Subset\Omega$ with $R<1.$ Then for any $0<\rho<\frac{R}{2}$ and $d>0$ we get that 
    \begin{equation*}
        \begin{aligned}
            \int_{B_{\rho}}\int_{B_{\rho}}|\log(\frac{u(x)+d}{u(y)+d})| \frac{dy dx}{|x+y|^N}\leq c \widetilde{M}^{2}(\rho^{N}&+\rho^{N+sp}d^{1-p}\int_{\RR^N\setminus B_R}\frac{u_{-}^{p-1}(y)+u_{-}^{q-1}(y)}{|y-x_0|^{N+sp}}dy
            \\
            &+\rho^{N+sq}d^{1-q}\int_{\RR^N\setminus B_R}\frac{u_{-}^{q-1}(y)}{|y-x_0|^{N+sq}}dy)
        \end{aligned}
    \end{equation*}
    for some $c(N,s,p,q,\mu_0)>0$ and $\widetilde{M}\equiv\widetilde{M}(\Omega^\prime)=1+(\|u\|_{L^{\infty}(\Omega^\prime)}+d)^{q-p}.$
\end{lemma}
\begin{sketch}
    For $\tau\geq 0$ set
    \begin{equation*}
        \begin{aligned}
            &\widetilde{H}(x,y,\tau)=G(\tau)=\frac{\tau^p}{\rho^{sp}}+\mu \frac{\tau^q}{\rho^{sq}},
            \\
            &\widetilde{h}(x,y,\tau)=g(\tau)=\frac{\tau^{p-1}}{\rho^{sp}}+\mu \frac{\tau^{q-1}}{\rho^{sq}}.
        \end{aligned}
    \end{equation*}
    Let $\phi\in C_{0}^{\infty}(B_{3\rho/2})$ satisfying $0\leq \phi\leq 1$ and $\phi=1$ in $B_\rho$ and $|D\phi|\leq 4/\rho.$  Use $\varphi(x)=\frac{\phi(x)^q}{g(u(x)+d)}$ as test function in $(Q_\mu)$ to get that 
    \begin{equation*}
        \begin{aligned}
            0&\leq \int_{\Omega}\varphi\leq I_1+I_2 \text{ where }
            \\
            I_{1}&=\int_{B_{2\rho}} \int_{B_{2\rho}}\frac{[u(x)-u(y)]^{p-1}}{|x-y|^{N+sp}}(\varphi(x)-\varphi(y))+\mu \frac{[u(x)-u(y)]^{p-1}}{|x-y|^{N+sp}}(\varphi(x)-\varphi(y)),
            \\
            I_{2}&=2\int_{\RR^N \setminus B_{ 2\rho}} \int_{B_{2\rho}}\frac{[u(x)-u(y)]^{p-1}}{|x-y|^{N+sp}}\varphi(x)+\mu \frac{[u(x)-u(y)]^{p-1}}{|x-y|^{N+sp}}\varphi(x).
        \end{aligned}
    \end{equation*}
    Next we put $a(x,y)=a_2=\mu$ and $[a]_{\alpha}=\|a\|_{L^\infty}=\mu_0$ in the argument \cite[Lemma 5.1]{BOS22}.  Observe that with this choice, \cite[(5.3)]{BOS22} is vacuously true. i.e.,  $G( u (y) + d )\leq (1+8^\alpha\mu_0)\;\widetilde{M}\; \widetilde{H}(x,y,u (y) + d).$ Then we argue similar to \cite[Lemma 5.1]{BOS22} and conclude the result.
\end{sketch}
Now similar to \cite[Corollary 5.2]{BOS22} we get the following result:
\begin{corollary}\label{logcor}
    Under the same assumptions of Lemma \ref{loglem} and let $d,\zeta>0,\xi>1$ and define 
    \begin{equation*}
        \begin{aligned}
            v:=\min\{(\log\frac{\zeta+d}{u+d})_+,\log \xi\},
        \end{aligned}
    \end{equation*}
    then we have 
    \begin{equation*}
        \begin{aligned}
            \intavg_{B_\rho} |v-(v)_{B_\rho}| dx \leq c \widetilde{M}^{2}(1&+\rho^{sp}d^{1-p}\int_{\RR^N\setminus B_R}\frac{u_{-}^{p-1}(y)+u_{-}^{q-1}(y)}{|y-x_0|^{N+sp}}dy
            \\
            &+\rho^{sq}d^{1-q}\int_{\RR^N\setminus B_R}\frac{u_{-}^{q-1}(y)}{|y-x_0|^{N+sq}}dy)
        \end{aligned}
    \end{equation*}
   for some $c(N,s,p,q,\mu_0)>0$ and $\widetilde{M}\equiv\widetilde{M}(\Omega^\prime)=1+(\|u\|_{L^{\infty}(\Omega^\prime)}+d)^{q-p}.$
\end{corollary}
First we recall a classical lemma  \cite[Lemma 7.1]{GIU03}.
\begin{lemma}
    \label{convlem}
   Let $\alpha>0$ and let $\{x_i\}$ be a sequence of positive real numbers such that 
 $ x_{i+1}\leq C B^{i}x_{i}^{1+\alpha}$ with $C>0$ and $B>1.$ If $x_0\leq C^\frac{-1}{\alpha}B^\frac{-1}{\alpha^2}$ we have $x_i\leq B^{\frac{-i}{\alpha}}x_0$ and hence in particular $\lim_{i\rightarrow\infty}x_i=0.$
\end{lemma}
Now we fix a ball $B_{2r}(x_0)\subset\Omega^{\prime}\Subset\Omega.$ Let $M=1+\|u\|_{L^\infty(\Omega^\prime)}^{q-p}$ and $\sigma\in(0,1/4]$ be a constant depending only on $N,s,p,q,\mu_0$ and $\|u\|_{L^{\infty}(\Omega^\prime)}$ that satisfies 
\begin{equation*}
    \begin{aligned}
        \sigma\leq \min\{\frac{1}{4},2^\frac{-2}{sp},6^{-\frac{4(q-1)}{sq}},\exp(-\frac{c_*M^3}{\nu_*})\}
    \end{aligned}
\end{equation*}
where the large constant $c_*(N,s,p,q,\mu_0)>0$ and the small constant $\nu_*(N,s,p,q,\mu_0,\|u\|_{L^{\infty}(\Omega^\prime)})>0$ are to be determined later. Then choose $0<\gamma(N,s,p,q,\mu_0,\|u\|_{L^{\infty}(\Omega^\prime)})<1$  satisfying 
\begin{equation*}
    \begin{aligned}
        \gamma\leq \min\{\log_\sigma(1/2),\frac{sp}{2(p-1)},\frac{sq}{2(q-1)},\log_\sigma(1-\sigma^{\frac{sq}{2(q-1)}})\}.
    \end{aligned}
\end{equation*}
We define 
\begin{equation*}
    \begin{aligned}
        \frac{K_{0}}{2}:=\sup_{B_r}|u|+[r^{sp}\int_{\RR^N\setminus B_r}\frac{|u(x)|^{p-1}+ |u(x)|^{q-1}}{|x-x_0|^{N+sp}}]^\frac{1}{p-1}+[r^{sq}\int_{\RR^N\setminus B_r}\frac{ |u(x)|^{q-1}}{|x-x_0|^{N+sq}}]^\frac{1}{q-1}>0
    \end{aligned}
\end{equation*}
and for $j\in\mathbb{N}\cup\{0\}$ we write
\begin{equation*}
    \begin{aligned}
        r_j=\sigma^j r, \text{ }B_j=B_{r_j}(x_0)\text{ and }K_j=\sigma^{\gamma j}K_0.
    \end{aligned}
\end{equation*}
Next we prove the oscillation lemma.
\begin{lemma}
    Let $u$ be a weak solution to the problem $(Q_\mu).$ Then we have for $j\in\mathbb{N}\cup\{0\}$
    \begin{equation*}
        \begin{aligned}
            \omega(r_j)=\osc_{B_j} u \leq K_j.
        \end{aligned}
    \end{equation*}
\end{lemma}
\begin{sketch}
    The proof goes by induction on $j.$ For $j=0$ it is obvious from the definition of $K_0.$ Now we assume the induction hypothesis is true for $i=0,\ldots,j$ for some fixed $j>0$ and want to show that it is true for $j+1.$ Without loss of generality we assume $\omega(r_{j+1})\geq \frac{K_{j+1}}{2}.$ Set
    \begin{equation*}u_j=\left\{
        \begin{aligned}
          &  u-\inf_{B_j} u \text{ \;\;\;\;if } \frac{|2 B_{j+1}\cap\{u\geq\inf_{B_j}(u+\frac{\omega(r_j)}{2}) \}|}{|2B_{j+1}|}\geq \frac{1}{2},\\
               &  \sup_{B_j} u-u \text{ \;\;\;\;if } \frac{|2 B_{j+1}\cap\{u\leq\inf_{B_j}(u+\frac{\omega(r_j)}{2}) \}|}{|2B_{j+1}|}\geq \frac{1}{2}.
        \end{aligned}\right.
    \end{equation*}
    Similar to \cite[step 2 ( Tail estimates) of Lemma 5.3]{BOS22} we get that 
    \begin{equation*}
        \begin{aligned}
            &r_j^{sp}\int_{\RR^N\setminus B_j}\frac{|u_j(x)|^{p-1}+|u_j(x)|^{q-1}}{|x-x_0|^{N+sp}} dx \leq c M \sigma^{-\gamma(p-1)} K_j^{p-1} \text{ and }
            \\
             &r_j^{sq}\int_{\RR^N\setminus B_j}\frac{|u_j(x)|^{q-1}}{|x-x_0|^{N+sp}} dx \leq c  \sigma^{-\gamma(q-1)} K_j^{q-1}
        \end{aligned}
    \end{equation*}
    for a constant $c(N,s,p,q)>0.$ We next apply Corollary \ref{logcor} to 
    \begin{equation*}
        \begin{aligned}
            v:=\min\{[\log(\frac{\omega(r_j)/2+d_j}{u_j+d_j})]_{+},k\}
        \end{aligned}
    \end{equation*}
    where $k=\log(\frac{1/2+\varepsilon}{3\varepsilon})>0$ and $d_j=\varepsilon K_j$ with $\varepsilon=\sigma^\frac{sq}{2(q-1)}.$ Thus similar to \cite[ step 3 of Lemma 5.4]{BOS22} we get that 
    \begin{equation*}
        \begin{aligned}
            \frac{|2 B_{j+1}\cap\{u_j\leq d_j \}|}{|2B_{j+1}|}\leq \frac{cM^{3}}{k} \leq \frac{c_*M^3}{\log(1/\sigma)}
        \end{aligned}
    \end{equation*}
    where $c$ and $c_*$ only depend on $N,s,p,q,\Omega$ and $\mu_0.$ Now we proceed with an iteration argument similar to \cite[step 4 of Lemma 5.4]{BOS22}. For $i=0,1,2,\ldots$ and for fixed $j$ we define
    \begin{equation*}
        \begin{aligned}
            \rho_i=(1+2^{-i})r_{j+1},\text{ }\tilde{\rho}_i=\frac{\rho_i+\rho_{i+1}}{2},\text{ }B^{i}=B_{\rho_i},\text{ }\tilde{B}_i=\tilde{B}_{\rho_i},
        \end{aligned}
    \end{equation*}
    and choose cut-off functions satisfying
    \begin{equation*}
        \begin{aligned}
            \phi_i\in C_0^\infty(\tilde{B}_i),\text{ }0\leq \phi_i\leq 1,\text{ }\phi_i\equiv 1\text{ on } B^{i+1} \text{ and } |D\phi_i|\leq \frac{2^{i+2}}{r_{j+1}}.
        \end{aligned}
    \end{equation*}
    Set 
    \begin{equation*}
        \begin{aligned}
         k_{i}=   (1+2^{-i})d_{j},\text{ }w_i=(k_i-u_j)_+,\text{ and  }A_i=\frac{|B^i\cap\{u_j<k_i\}|}{|B^i|}.
        \end{aligned}
    \end{equation*}
    Put $a_{1,j}=a_{2,j}=\mu$ in the arguments of \cite[step 4 of Lemma 5.4]{BOS22} and define 
    \begin{equation*}
        \begin{aligned}
            G(\tau)=\frac{\tau^p}{r_{j+1}^{sp}}+\mu \frac{\tau^q}{r_{j+1}^{sq}}.
        \end{aligned}
    \end{equation*}
    Then for $\kappa=\min\{\frac{p_s^*}{p},\frac{q_s^*}{q}\}>1,$ we observe that for any function $f\in \mathcal{W} (B_{r_{j+1}})\cap L^{\infty}(B_{r_{j+1}})$ we get that 
    \begin{equation*}
        \begin{aligned}
          [\intavg_{B_{r_{j+1}}} ( G(f))^{\kappa}]^{\frac{1}{\kappa}}\leq (1+\|f\|^{q-p}_{L^{\infty}(B_{r_{j+1}})}) [\intavg_{B_{r_{j+1}}} ( G(f))^{\kappa}]^{\frac{1}{\kappa}}.
        \end{aligned}
    \end{equation*}
    Then similar to \cite[Lemma 2.5]{BOS22} we get that 
    \begin{equation*}
        \begin{aligned}
            &\left[\intavg_{B_{r_{j+1}}} ( G(f))^{\kappa}\right]^{\frac{1}{\kappa}}
            \\
            &\leq c (1+\|f\|^{q-p}_{L^{\infty}(B_{r_{j+1}})})\left[\intavg_{B_{r_{j+1}}}\int_{B_{r_{j+1}}}H(x,y,|f(x)-f(y)|)\frac{dxdy}{|x-y|^N}
            +\intavg_{B_{r_{j+1}}} G(f)\right]
        \end{aligned}
    \end{equation*}
    for a constant $c(N,s,p,q,\Omega)>0.$ Then similar to \cite[(5.30))]{BOS22} we get that 
    \begin{equation*}
        \begin{aligned}
            A_{i+1}^{1/\kappa}G(k_i-k_{i+1})\leq cM \intavg_{B^{i+1}}\int_{B^{i+1}} H(x,y,|w_i(x)-w_i(y)|)\frac{dxdy}{|x-y|^N} + cM G(d_j) A_i.
        \end{aligned}
    \end{equation*}
    Similar to \cite[(5.31)]{BOS22} using lemma \ref{caclem} and the tail estimates we get that 
    \begin{equation*}
        \begin{aligned}
         & \intavg_{B^{i+1}}\int_{B^{i+1}} H(x,y,|w_i(x)-w_i(y)|)\frac{dxdy}{|x-y|^N}   
         \\
         &\leq  c\intavg_{B_{i}} \int_{B_{i}} H(x,y,|(w_i(x)+w_i(y))(\phi_i(x)-\phi_i(y))|)\frac{dx\;dy}{|x-y|^{N}}
           \\
           & + \left(\sup_{x\in\tilde{B}_i}\int_{\RR^{N}\setminus B_{i}}h(x,y,w_i(y))\frac{dy}{|x-y|^{N}})\right)\intavg_{B_{i}}w_i(x)\phi_i^{q}(x)\thinspace dx + \intavg_{B_{i}} w_i(x)\phi_i^{q}(x)\thinspace dx.
        \end{aligned}
    \end{equation*}
    From \cite[(5.32),(5.33) and (5.34)]{BOS22} we get that 
    \begin{equation*}
        \begin{aligned}
            &c\int_{B_{i}} \int_{B_{i}} H(x,y,|(w_i(x)+w_i(y))(\phi_i(x)-\phi_i(y))|)\frac{dx\;dy}{|x-y|^{N}}\leq c 2^{iq} G(d_j) A_i
            \\
            & \intavg_{B_{i}} w_i(x)\phi_i^{q}(x)\thinspace dx. \leq c d_j A_i
            \\
            & \left(\sup_{x\in\tilde{B}_i}\int_{\RR^{N}\setminus B_{i}}h(x,y,w_i(y))\frac{dy}{|x-y|^{N}})\right) \leq c 2^{i(N+sq)}M\frac{G(d_j)}{d_j}.
        \end{aligned}
    \end{equation*}
    Observe that $d_j\leq G(d_j)$ as $\gamma\leq  \frac{sp}{2(p-1)}$ by definition. Thus we have 
    \begin{equation*}
        \begin{aligned}
            \intavg_{B_{i}} w_i(x)\phi_i^{q}(x)\thinspace dx. \leq c G(d_j) A_i.
        \end{aligned}
    \end{equation*}
    So a similar calculation yields
    \begin{equation*}
        \begin{aligned}
            A_{i+1}\leq c_0 2^{i\kappa(N+sq+2q)} M^{2\kappa} A_i.
        \end{aligned}
    \end{equation*}
    Set 
    \begin{equation*}
        \begin{aligned}
            \nu_*=(c_0M^{-2\kappa})^{\frac{-1}{\kappa-1}}2^\frac{-\kappa(N+sq+2q)}{(\kappa-1)^2}
        \end{aligned}
    \end{equation*}
    to get that 
    \begin{equation*}
        \begin{aligned}
       A_0=     \frac{|2B_{j+1}\cap\{u_j\leq 2 d_j\}|}{|2B_{j+1}|} \leq \frac{c_* M^3}{\log(1/\sigma)}\leq \nu_*
        \end{aligned}
    \end{equation*}
which holds by the definition of $\sigma.$ By lemma \ref{convlem}, this implies $A_i\leq 2^{\frac{-i}{\kappa-1}}A_0\leq 2^{\frac{-i}{\kappa-1}}$ and consequently we get that $\lim_{i\rightarrow\infty} A_i=0$ independent of $\mu$ which means $u_j\geq d_j$ a.e. in $B_{j+1}.$ Since $\|u\|_{L^{\infty}(\Omega)}$ is independent of $\mu$   we can conclude that 
\begin{equation*}
    \begin{aligned}
        \omega(r_{j+1})\leq (1-\varepsilon)K_j\leq (1-\sigma^{\frac{sq}{2(q-1)}}) \sigma^{-\gamma} K_{j+1}\leq K_{j+1}.
    \end{aligned}
\end{equation*}
This completes the proof.
\end{sketch}
Combining all the arguments we get the following theorem:
\begin{theorem}\label{careghld}
Under the same assumptions as Theorem \ref{careg} when $f(x)\equiv 1$, and for \( s \in \left(0, \frac{1}{q} \right) \), there exists \(\alpha_0 \in (0,1)\) such that \(\zmu \in C_{\mathrm{loc}}^{0,\alpha_0}\) and for any $K\Subset \Omega$
\begin{equation*}
   \relax \|\zmu\|_{C^{\alpha_{0}}(K)} \leqslant C(K, N,s,p,q,\sigma,\mu_0)
\end{equation*}
independent of $\mu.$  
\end{theorem}
Next in this appendix, we present the superposition principle, Lewy-Stampachia inequality, and a technical lemma. These propositions, already established for the fractional $p$-Laplacian in \cite[Preliminaries]{IMS20}, can be extended to the fractional $(p,q)$-Laplacian framework. The variational approach of the proof  relies on the operator's monotonicity, without any crucial dependence on scaling and homogeneity properties. Owing to the straightforward adaptability afforded by this approach, we refrain from providing detailed proofs. We commence our discussion with the superposition principle, a pivotal element in our analytical framework.

\begin{proposition}(Superposition Principle)\label{superpostn}
    Let $\Om$ be bounded, $u\in \tw(\Om)$, $v\in L_{\loc}^{1}(\RR^N),$ $V=\supp(u-v)$ satisfy
    \begin{itemize}
        \item[i)] $\Om \subset\subset \RR^N\setminus V;$
        \item[ii)] $\int_{V}\frac{|v(x)|^{t-1}}{(1+|x|)^{N+ts}} \thinspace dx < \infty$ \text{ for } t=p,q.
    \end{itemize}
    Set for all $x \in \RR^N,$
    \begin{equation*}
        w(x)=\left\{\begin{aligned}
            &u(x) \text{ if } x \in V^{c}
            \\
            & v(x) \text{ if } x \in V.
        \end{aligned}\right.
    \end{equation*}
    Then  $w\in \tw(\Om)$ and satisfies in $\Om$
   \begin{equation*}
       \begin{aligned}
           \fp w(x) + \fq w(x) = &\fp u(x)+ \fq u (x) 
           \\&+\sum_{t=p,q} 2 \int_{\bR} \frac{(u(x)-v(y))^{t-1}-(u(x)-u(y))^{t-1}}{|x-y|^{N+ts}}. 
       \end{aligned}
   \end{equation*}
\end{proposition}

Next we shall discuss a generalisation of the Lewy-Stampacchia type inequality for fractional $(p,q)$ Laplace operator. Our proof is motivated by the insights in \cite[Page no 12]{IMS20} and the abstract Lewy Stampacchia result \cite{GM15}.  We introduce a partial ordering on the dual space $\mathcal{W}^{\prime}(\Om)$ by defining the positive cone
\begin{equation*}
	\mathcal{W}^{\prime}(\Om)_+=\lbrace L\in \mathcal{W}^{\prime}(\Om)\thinspace : \langle L {,} \varphi
	  \rangle \ge 0 \thickspace \text{for all}\thickspace \varphi \in \mathcal{W}_{0}(\Om)_{+}  \rbrace .
\end{equation*}
Using the ideas in \cite{FSV15}, we can prove that $C_{c}^{\infty}(\Om)$ is dense in $\mathcal{W}_0(\Omega)$. Now using Riesz representation theorem  any $L \in \mathcal{W}^{\prime}(\Om)_+ $ can be represented as a positive Radon measure on $\Om$. Then, the order dual is defined as:
\begin{equation*}
	\mathcal{W}^{\prime}_{\leqslant}(\Omega)=\{L_{1}-L_{2}\thickspace : \thickspace L_{1},L_{2}\in  \mathcal{W}'(\Om)_{+} \}.
	\end{equation*}
We define energy functional $E:\mathcal{W}_0(\Omega)\rightarrow \mathbb{R}$ as
\begin{equation*}
	E(u)=\frac{1}{p}[u]^{p}_{s,p}+\frac{1}{q}[u]^{q}_{s,q}.
\end{equation*}
	\begin{lemma}\label{LSlem} (Lewy-Stampacchia). Let $\Omega \subseteq \mathbb{R}^N$ be bounded, $\varphi, \psi \in 	\mathcal{W}_{\mathrm{loc}}\left(\mathbb{R}^N\right)$ be such that
		\begin{itemize}
			\item[i)] $\fp \varphi+\fq \varphi,\fp \psi+ \fq \psi \in \mathcal{W}^{\prime}_{\leqslant}(\Omega),$
			\item[ii)] $[\varphi, \psi]:=\left\{v \in \mathcal{W}_0(\Omega): \varphi \leqslant v \leqslant \psi\right\} \neq \emptyset,$
		\end{itemize} 
		Then there exists a unique solution $u \in 	\mathcal{W}_0(\Omega)$ to the problem
		\begin{equation*}
			\begin{aligned}
				\min _{v \in[\varphi, \psi]} E(v),
			\end{aligned}
		\end{equation*}
		and it satisfies
		\begin{equation*}
			\begin{aligned}
				0 \wedge \left(\fp\psi +\fq \psi\right) \leq \fp u +\fq u \leq 0 \vee \left(\fp\va +\fq \va\right)\text { in } \Omega.
			\end{aligned}
		\end{equation*}
		
	\end{lemma}

\begin{proof}
    Observe that $E$ is convex and coercive in $\wo(\Omega).$ Repeating the discussion of \cite[Page 266]{GM15}, we can show that $E$ is sub-modular. i.e. 
    \begin{equation*}
        \begin{aligned}
            E(u \vee v) + E(u \wedge v) \leqslant E(u)+ E(v) \;\;\; \forall u,v \in \wo(\Om).
        \end{aligned}
    \end{equation*}
    The strict convexity and coercivity of the functional $E$ ensure the existence of the unique solution to the minimization problem.  Further, the submodularity and strict convexity of $E$ imply that its differential, $\fp + \fq$, is a strictly $\mathcal{T}-$monotone map i.e. 
    \begin{equation*}
        \begin{aligned}
            \langle \fp u + \fq u - \fp v - \fq v, (u-v)_{+}\rangle >0 \text{ unless } v\leqslant u.  \end{aligned}
    \end{equation*}
 Now the proof follows a similar way as presented in \cite[Lemma 2.5]{IMS20}. 
\end{proof}

The following proposition establishes a crucial estimate for functions that exhibit local boundedness by a suitable multiple of $\doos$. The significance of this estimate and the necessity of the employed decomposition are elaborated in  \cite[Remark 2.8]{IMS20}. Leveraging the monotonicity of our operator and noting that the proof of  \cite[Proposition 2.7]{IMS20} doesn't hinge on scaling or homogeneity properties, we establish this estimate analogously.
\begin{proposition}\label{maxminlem}
    Let $\Om$ be bounded, and $u \in \two(D_R)$ satisfy $\fp u + \fq u \in \mathcal{W}^{\prime}_{\leqslant}(D_R).$
    \begin{itemize}
        \item[i)] Suppose $u\geq m \doos$ in $D_{2R}.$ There exist $C_{2,t}(N,t,s)>0$ and $C_{\varepsilon,t}(N,t,s,\varepsilon)>0$ for all $\varepsilon>0$ such that in $D_R$
        \begin{eqnarray*}
            \begin{array}{rl}
               \fp (u \vee m\doos) +  \fq (u \vee m\doos) 
                \geq \fp u &+ \fq u 
               +\sum_{t=p,q}\Big( \frac{-\varepsilon}{R^s} \|\frac{u}{\doos}-m\|_{L^{\infty}(D_R)}^{t-1}
               \\
               &- C_{\varepsilon,t} \tail_{t-1} ((m-\frac{u}{\doos})_{+},2R)^{t-1}
               \\
               &-C_{2,t} |m|^{t-2} \tail_1((m-\frac{u}{\doos})_{+},2R) \Big).
            \end{array}
        \end{eqnarray*}

        \item[ii)] Suppose $u\leq M \doos$ in $D_{2R}.$ There exist $C'_{2,t}(N,t,s)>0$ and  $C'_{\varepsilon,t}(N,t,s,\varepsilon)>0$  for all $\varepsilon>0$ in $D_R$
        \begin{eqnarray*}
            \begin{array}{rl}
                \fp (u \wedge M\doos) +  \fq (u \wedge M\doos) 
                \leq \fp u &+ \fq u 
               + \sum_{t=p,q}\Big( \frac{\varepsilon}{R^s} \|M-\frac{u}{\doos}\|_{L^{\infty}(D_R)}^{t-1}
               \\
               &+ C'_{\varepsilon,t} \tail_{t-1} ((\frac{u}{\doos}-M)_{+},2R)^{t-1}
               \\
               &+C'_{2,t} |M|^{t-2} \tail_1((\frac{u}{\doos}-M)_{+},2R) \Big).
            \end{array}
        \end{eqnarray*}

    \end{itemize}
\end{proposition}

\bibliographystyle{abbrv}
	\bibliography{ref}

\end{document}